\documentclass[12pt,a4paper]{article}
\usepackage{amssymb,latexsym,bm,amsmath,amsthm}
\usepackage[dvips]{color}
\usepackage{array}
\usepackage[section]{placeins}
\usepackage{arcs}
\usepackage{epsfig}
\graphicspath{{eps_figs/}}

\newtheorem{thm}{Theorem}[section]
\newtheorem{prop}[thm]{Proposition}
\newtheorem{defn}[thm]{Definition}
\newtheorem{lem}[thm]{Lemma}
\newtheorem{rmk}[thm]{Remark}
\newtheorem{cor}[thm]{Corollary}

\newcommand{\argmax}{\mathop{\rm argmax}}

 \def\O{{\cal O}}
\def\tr{{\rm trace}}
\def\rmspan{{\rm span}}

\def\0{{\bm 0}}\def\vec{{\sf vec}}
\def\cov{{\rm Cov}} \def\acov{{\rm aCov}}

\def\ptilde{{\tilde p}}\def\qtilde{{\tilde q}}

\def\Ahat{{\widehat A}}\def\Bhat{{\widehat B}}

\definecolor{blue}{rgb}{0.1,0.1,1}
\definecolor{darkblue}{rgb}{0,0,0.5}

\begin{document}
\setcounter{page}{1}
\def\thepage{0-\arabic{page}}
\setcounter{page}{1} \setcounter{equation}{0}
\def\thepage{\arabic{page}}

\title{On Multilinear Principal Component Analysis of Order-Two Tensors}
\author{Hung Hung$^{a}$, Pei-Shien Wu$^{b}$, I-Ping Tu$^{c,}$\footnote{Corresponding author.
\textit{Email address}: iping@stat.sinica.edu.tw }, and Su-Yun Huang$^c$\\
\small $^a$Institute of Epidemiology and Preventive Medicine,
National Taiwan University\\
\small $^b$Department of Mathematics, National Taiwan University\\
\small $^c$Institute of Statistical Science, Academia Sinica}
\date{\empty} \maketitle

\begin{abstract}
Principal Component Analysis (PCA) is a commonly used tool for
dimension reduction in analyzing high dimensional data; Multilinear
Principal Component Analysis (MPCA) has the potential to serve the
similar function for analyzing tensor structure data. MPCA and other
tensor decomposition methods have been proved effective to reduce
the dimensions for both real data analyses and simulation studies
(Ye, 2005; Lu, Plataniotis and Venetsanopoulos, 2008; Kolda and
Bader, 2009; Li, Kim and Altman, 2010). In this paper, we
investigate MPCA's statistical properties and provide explanations
for its advantages. Conventional PCA, vectorizing the tensor data,
may lead to inefficient and unstable prediction due to its extremely
large dimensionality. On the other hand, MPCA, trying to preserve
the data structure, searches for low-dimensional multilinear
projections and decreases the dimensionality efficiently. The
asymptotic theories for order-two MPCA, including asymptotic
distributions for principal components, associated projections and
the explained variance, are developed. Finally, MPCA is shown to
improve conventional PCA on analyzing the {\sf Olivetti Faces} data
set, by constructing more module oriented basis in reconstructing
the test faces.
\end{abstract}


\section{Introduction}
Dimension reduction is a key step for high dimensional data
analysis. Principal component analysis (PCA) is probably the most
commonly used method for dimension reduction. Given $n$ observations
on $m$ variables, PCA calculates the $m\times m$ covariance matrix
and solves the  eigenvalue decomposition problem for the covariance
matrix. The goal is to choose a smaller set  of eigenvectors as a
new coordinate system so that the newly transformed variables can
retain the most data variation. This PCA approach has been widely
applied in many scientific fields for dimension reduction and
compact data representation (Jolliffe, 2002), where the collected
data are organized in an $n\times m$ design matrix with each row
representing an observation and each column a variable.

When data are tensor objects, traditional analysis vectorizes each
of the tensor objects into a long vector and arranges these
vectorized objects in a design matrix form. Subsequent analysis is
followed in the usual way. Nevertheless, this approach usually
produces a large number of variables, where the available sample
size is relatively small, and many existing statistical methods fail
to apply. For a typical example, like the {\sf Olivetti Faces} data
set to be used in an experimental study later, there are 400 images
each with $64\times 64$ pixels. Vectorizing each image leads to a
design matrix of size $400\times 4096$, which the variable dimension
$m$ largely exceeds the sample size $n$.

One strategy to overcome this difficulty is to take advantage of the
natural tensor structure of the data. Singular value decomposition
(SVD) is an example. Given a $p\times q$ matrix $X$ which can be
treated as an order-two tensor, SVD can decompose two directional
spaces simultaneously: $X=USV^T=\sum_{i=1}^{p\wedge q}s_iu_i v_i^T$,
where $U=[u_1,..., u_p]\in \mathfrak{R}^{p\times p}$ and
$V=[v_1,...,v_q]\in \mathfrak{R}^{q\times q}$ are, respectively, the
left and right singular vectors, $S$ is a diagonal matrix of size
$p\times q$ with diagonal elements $\{s_1,...,s_{(p\wedge q)}\}$.
The dimension can be reduced when the index $i$ is properly
truncated. De Lathauwer, De~Moor and Vandewalle (2000a) then
generalized the SVD to high-order SVD (HOSVD) for a given $N^{\rm
th}$-order tensor object $A\in \Re^{I_1\times\cdots\times I_N}$.
Further, they formulated the problem of ``best
rank-$(R_1,\cdots,R_N)$ approximation of higher-order tensors'' in
the least-squares sense, and discussed many algorithms to achieve
this task (De Lathauwer, De~Moor and Vandewalle, 2000b).

Later, Yang et al. (2004) proposed two-dimensional PCA (2DPCA) for
analyzing image data, which are order-two tensors. An improved
two-directional two-dimensional PCA (${\rm (2D)}^2$PCA) was
developed in Zhang and Zhou (2005), which was shown to perform
better than 2DPCA through simulation studies. Ye (2005) formulated
the problem of generalized low rank approximation of matrices,
which can be treated as a sample extension of the best
rank-$(R_1,R_2)$ approximation for order-two tensors in
De~Lathauwer, De~Moor and Vandewalle (2000b). Lu, Plataniotis and
Venetsanopoulos~(2008) further generalized the work of Ye (2005)
and proposed multilinear PCA (MPCA) for tensor objects of
arbitrary orders. There are other tensor decomposition methods for
dimension reduction. For instance, Kolda and Bader (2009) provided
a general overview of current development of tensor decomposition
methods for unsupervised learning, their applications, and
available softwares; Li, Kim and Altman. (2010) considered the
tensor decomposition methods for supervised learning such as
regression and classification.

Similar to conventional PCA, the goal of MPCA is to look for
low-dimensional multilinear projection for tensor objects that
captures the most data variation. Back to the example of {\sf
Olivetti Faces}, one eigenvector in conventional PCA creates an
image basis element that contains 4095 free parameters. By contrast,
one image basis element in MPCA or (2D)$^2$PCA, which involves the
Kronecker product of a column vector and a row vector, contains 126
free parameters. From the viewpoint of the number of parameters
required to specify one basis element, MPCA is expected to perform
better than conventional PCA, when the sample size is small to
moderate, like this {\sf Olivetti Faces} example. Compared to
(2D)$^2$PCA, MPCA has the advantage of capturing more data variation
by the chosen image basis, because of its specific criterion. MPCA
has been successfully applied in real data analysis and checked by
simulations (Ye, 2005; Lu et al., 2008). Yet, to our best knowledge,
there is neither statistical justification nor asymptotic study for
MPCA.

In this paper, we try to establish some relevant properties of
order-two MPCA from a statistical point of view. Our study is based
on the following model:
\begin{equation}
X=\mu+ A_0UB_0^T+ \varepsilon,\label{model}
\end{equation}
where $\mu\in \Re^{p\times q}$ is the mean parameter of $X$, $A_0\in
\Re^{p\times p_0}$ and $B_0 \in \Re^{q\times q_0}$ with $p_0\leq p$
and $q_0\leq q$ are non-random basis matrices, $U\in\Re^{p_0\times
q_0}$ is a random coordinate matrix with $E[U]=0$ and a strictly
positive definite covariance matrix
$\cov(\vec(U))=T\in\Re^{m_0\times m_0}$, where $m_0=p_0q_0$ and
$\vec(\cdot)$ is the operator that stacks the columns of a tensor
into a long vector. The error term $\varepsilon\in\Re^{p\times q}$
is a random matrix independent of $U$ and with $E[\varepsilon]=0$
and $\cov(\vec(\varepsilon))=\sigma^2I_m$, where $m=pq$. Under
model~(\ref{model}) which characterizes the tensor structure of $X$,
we justify the validity of MPCA. Asymptotic properties of MPCA are
rigorously developed, including asymptotic distributions for
principal components, associated projections and the explained
variance. It is also shown that MPCA is asymptotically more
efficient than (2D)$^2$PCA in estimating the target dimension
reduction subspace. Furthermore, a test of dimensionality is
developed, based on the derived asymptotic results.

This paper is organized as follows.  Section 2 presents some
properties of the estimation for the target subspace and a test for
its dimensionality. The relations between MPCA and both conventional
PCA and (2D)$^2$PCA are also discussed in this section. In
section~3, the asymptotic theory of MPCA is developed. In section 4,
the performance of MPCA and its comparison with conventional PCA is
demonstrated by analyzing the {\sf Olivetti Faces} data set. The
paper ends with a brief discussion. Technical proofs of main results
are deferred to the Appendix.

\section{MPCA}
MPCA, as a dimension reduction algorithm, is originally designed to
search basis matrices $\{A,B\}$ and coordinate matrices $U_i$'s that
best approximate the observed data $X_i$ as $AU_iB^T$ for
$i=1,\dots,n$. Although many simulation studies and real data
analyses in literature support the usage of MPCA and multilinear
tensor decomposition (Ye, 2005; Lu et al., 2008; Kolda and Bader,
2009; Li et al., 2010), there is no theoretical study from the
statistical point of view. { Let $\otimes$ be the Kronecker product.
Then, there is an equivalent formula for
 model~(\ref{model})
\begin{equation}
X=\mu+ A_0UB_0^T+ \varepsilon ~~\Leftrightarrow~~
\vec(X-\mu)=(B_0\otimes A_0)\vec(U)+\vec(\varepsilon)\label{model1}
\end{equation}
by the fact that $\vec(A_0UB_0^T)=(B_0\otimes A_0)\vec(U)$. Without
loss of generality, we may assume that $A_0$ and $B_0$ are
orthogonal matrices, i.e., $A_0^TA_0=I_{p_0}$ and
$B_0^TB_0=I_{q_0}$. Model~(\ref{model}) thus ensures that, without
considering the error term $\varepsilon$, the columns and rows of
$(X-\mu)$ belong to $\rmspan(A_0)$ and $\rmspan(B_0)$, respectively,
and $\vec(X-\mu)$ belongs to the subspace
$\rmspan(B_0)\otimes\rmspan(A_0)=\rmspan(B_0\otimes A_0)$. It is
then reasonable to estimate $\rmspan(B_0\otimes A_0)$ for follow-up
analysis such as data compression, pattern recognition, regression
analysis, etc.} In this section, we show that, under model
(\ref{model}), MPCA actually attempts to extract a basis pair
$\{A,B\}$ targeting the subspace $\rmspan(B_0\otimes A_0)$.
Proposition~\ref{lrapprox} below proves the existence of a solution
pair $\{A,B\}$. Proposition~\ref{relationship} summarizes that the
inclusion relation between $\rmspan(A)$ (resp., $\rmspan(B)$) and
the target dimension reduction subspace $\rmspan(A_0)$ (resp.,
$\rmspan(B_0)$), depends on the size comparison between the
specified dimensionality $\ptilde$ (resp., $\qtilde$) and $p_0$
(resp., $q_0$). Recognizing the important roles of $\ptilde$ and
$\qtilde$, we construct a hypothesis test for choosing $\ptilde$ and
$\qtilde$. These works justify the usage of MPCA in extracting the
relevant basis for subsequent analysis, provided the data has a
natural tensor structure.

\subsection{Estimation}

Let $\{X_i\}_{i=1}^n$ be the collected data set which are assumed to
be random copies of a random matrix $X\in \Re^{p\times q}$. MPCA
aims to extract the basis pair that best approximate
$\{X_i\}_{i=1}^n$ while preserving the tensor structure of them. In
particular, for a pre-specified dimensionality $(\ptilde,\qtilde)$,
Ye (2005) proposed a criterion to find
$A\in\O_{p,\ptilde},B\in\O_{q,\qtilde}$, and $\{U_i\}_{i=1}^n$ that
minimize
\begin{equation}\label{obj}
\frac{1}{n}\sum_{i=1}^n\|(X_i-\bar X)-AU_iB^T\|_F^2,
\end{equation}
where $\bar X=\frac1n\sum_{i=1}^n X_i$ is the sample mean matrix,
$\|\cdot\|_F$ is the Frobenius norm of a matrix, and
$\O_{\ell,\tilde\ell}$ is the collection of all orthogonal matrices
$M$ of size $\ell\times\tilde\ell$ such that $M^TM=I_{\tilde\ell}$.
Note that the objective function (\ref{obj}) can be expressed as
\begin{equation}\label{obj2}
\frac{1}{n}\sum_{i=1}^n\|\vec(X_i-{\bar X})-(B\otimes
A)\vec(U_i)\|^2_2.
\end{equation}
If we replace $(B\otimes A)$ by $\Gamma$ in (\ref{obj2}), the
minimization problem then becomes the conventional PCA. From this
viewpoint, MPCA can be treated as a constrained PCA with the
tensor constraint $\Gamma=(B\otimes A)$, where
$A\in\O_{p,\ptilde}$ and $B\in\O_{q,\qtilde}$. The following
theorem established in Ye (2005) characterizes some useful
properties of the solutions of the minimization problem
(\ref{obj}). In the rest of discussion, $P_M$ denotes the
orthogonal projection matrix onto $\rmspan(M)$ and $Q_M=I-P_M$.
\begin{thm}\label{thm.ye}(Ye, 2005)
Let $\widehat A,\widehat B,\{\widehat U_i\}_{i=1}^n$ constitutes a
minimizer for (\ref{obj}) under the dimensionality
$(\ptilde,\qtilde)$. Then,
\begin{itemize}
\item[(a)]
$\widehat U_i=\widehat A^T(X_i-{\bar X})\widehat B$.
\item[(b)]
$\{\widehat A, \widehat B\}$ is the maximizer of
$\frac{1}{n}\sum_{i=1}^n\|A^T(X_i-{\bar X})B\|_F^2$.
\item[(c)]
$\widehat A$ consists of the leading $\ptilde$ eigenvectors of
$\frac{1}{n}\sum_{i=1}^n(X_i-{\bar X})P_{\widehat B}(X_i-{\bar
X})^T$, and $\widehat B$ consists of the leading $\qtilde$
eigenvectors of $\frac{1}{n}\sum_{i=1}^n(X_i-{\bar
X})^TP_{\widehat A}(X_i-{\bar X})$.
\end{itemize}
\end{thm}
Similarly, we can define a population version of (\ref{obj}):
$E\|(X-\mu)-AUB^T\|_F^2$, and the corresponding minimizer should
follow Theorem \ref{thm.ye} such that the minimizer over
$A\in\O_{p,\ptilde}$ and $B\in\O_{q,\qtilde}$,  is equivalent to the
maximizer of the maximization problem:
\begin{equation}\label{max_prob}
\argmax_{A\in\O_{p,\ptilde},B\in\O_{q,\qtilde}} E\|A^T(X-\mu)B\|_F^2
=\argmax_{A\in\O_{p,\ptilde},B\in\O_{q,\qtilde}} \tr\left\{(B\otimes
A)^T\Sigma (B\otimes A)\right\},
\end{equation}
where $\Sigma=\cov(\vec(X))$. The following proposition gives the
existence of the solution.

\begin{prop}\label{lrapprox}
For a fixed but arbitrary positive semi-definite matrix $\Sigma$
of size $pq\times pq$, solution(s) to the maximization problem
(\ref{max_prob}) exists.
\end{prop}

Note that we do not need the model assumption~(\ref{model}) for
Proposition~\ref{lrapprox}. Also note that
Proposition~\ref{lrapprox} applies to problem (\ref{obj}) as well by
replacing $\Sigma$ with its sample estimate $S_n$, the sample
covariance matrix of $\{\vec(X_i)\}_{i=1}^n$, and by rephrasing the
maximization problem into the equivalent minimization problem. With
the existence of the maximizer in (\ref{max_prob}) we can formally
define the tensor principal components and the MPCA subspace.

\begin{defn}\label{def}
For a pre-specified dimensionality $(\ptilde,\qtilde)$, let
$\{A,B\}$ be the unique solution to the maximization
problem~(\ref{max_prob}), where $A$ and $B$ can be expressed in
their columns as $A=[a_1,\dots,a_\ptilde]$ and
$B=[b_1,\dots,b_\qtilde]$. We call $\{b_j\otimes a_i: 1\leq i\leq
\ptilde,1\leq j\leq \qtilde\}$ the tensor principal components, and
$\rmspan(B\otimes A)$ the MPCA subspace of dimensionality
$(\ptilde,\qtilde)$.
\end{defn}

Using similar arguments as in Theorem~\ref{thm.ye}~(c), we have that
$A$ and $B$ consist of the leading $\ptilde$ and $\qtilde$
eigenvectors of $E[(X-\mu)P_B(X-\mu)^T]$ and $E[(X-\mu)^T
P_A(X-\mu)]$, respectively. Since
\begin{eqnarray}
E[(X-\mu)P_B(X-\mu)^T]=\sum_{j=1}^{\qtilde}E[(X-\mu) (b_jb_j^T)
(X-\mu)^T] =\sum_{j=1}^{\qtilde}(b_j\otimes
I_p)^T\Sigma(b_j\otimes I_p),\label{def2.3_eq1}\\
E[(X-\mu)^TP_A(X-\mu)]=\sum_{i=1}^{\ptilde}E[(X-\mu)^T (a_ia_i^T)
(X-\mu)]=\sum_{i=1}^{\ptilde} (a_i\otimes I_p)^T\Sigma(a_i\otimes
I_p),\label{def2.3_eq2}
\end{eqnarray}
equivalently, $\{A,B\}$ consist of the leading solutions of the
system of stationary equations
\begin{eqnarray*}
\left(\sum_{j=1}^{\qtilde}(b_j\otimes I_p)^T\Sigma(b_j\otimes I_p)\right)a_i&=&\lambda_ia_i, ~i=1,\cdots,\ptilde,\\
\left(\sum_{i=1}^{\ptilde}(I_q\otimes a_i)^T\Sigma(I_q\otimes
a_i)\right)b_j&=&\xi_j b_j, ~j=1,\cdots,\qtilde,
\end{eqnarray*}
over $A\in \mathcal{O}_{p\times\ptilde}$ and $B\in
\mathcal{O}_{q\times\qtilde}$, where the ordering is determined by
the corresponding eigenvalues $\lambda_1\geq\cdots
\geq\lambda_{\ptilde}\geq 0$ and
$\xi_1\geq\cdots\geq\xi_{\qtilde}\geq 0$.

\begin{rmk}\label{rmk.def}
Obviously $\lambda_i$'s, $a_i$'s, $\xi_j$'s and $b_j$'s depend on
$\Sigma$. Besides such dependence, they also depend on the
dimensionality $(\ptilde,\qtilde)$. A more precise notation for
them should be $\lambda_i(\Sigma,\ptilde,\qtilde)$,
$a_i(\Sigma,\ptilde,\qtilde)$, $\xi_j(\Sigma,\ptilde,\qtilde)$ and
$b_j(\Sigma,\ptilde,\qtilde)$. However, for notation simplicity,
we use $\lambda_i$, $a_i$, $\xi_j$ and $b_j$, unless we want to
emphasize on their dependence on $(\Sigma,\ptilde,\qtilde)$.
\end{rmk}

From Remark~\ref{rmk.def}, for any fixed $(\ptilde,\qtilde)$, we
could define the sample analogues $\{\widehat A,\widehat B\}$,
$\hat\lambda_i$'s, and $\hat\xi_j$'s by replacing $\Sigma$ with the
sample covariance matrix $S_n$. In the rest of the discussion, with
pre-specified dimensionality $(\ptilde,\qtilde)$, we denote the
solution of (\ref{max_prob}) by $\{A,B\}$ and the population tensor
principal components by $(B\otimes A)$, and the corresponding sample
analogues by $\{\Ahat,\Bhat\}$ and $(\Bhat\otimes \Ahat)$. Finding
principal components in conventional PCA is equivalent to an
eigenvalue-problem. However, there is no explicit solution of
$\{\widehat A, \widehat B\}$ for MPCA; therefore, an algorithm was
proposed. The GLRAM algorithm of Ye (2005) to obtain $\{\widehat
A,\widehat B\}$ is summarized below.\\

\noindent \textbf{GLRAM (Ye, 2005):} Given a random initial
$A^{(0)}\in \mathcal{O}_{p\times \ptilde}$. For $k=1,2,\cdots,$
\begin{itemize}
\item[1.]
Obtain the maximizer $B^{(k+1)}=\arg\max_{B\in \mathcal{O}_{q\times
\qtilde}}\frac{1}{n}\sum_{i=1}^n\|A^{(k)T}(X_i-\bar X)B\|_F^2$.
\item[2.]
Obtain the maximizer $A^{(k+1)}=\arg\max_{A\in \mathcal{O}_{p\times
\ptilde}}\frac{1}{n}\sum_{i=1}^n\|A^{T}(X_i-\bar X)B^{(k+1)}\|_F^2$.
\item[3.]
Repeat Steps~1-2 until there is no significant difference between
$\frac{1}{n}\sum_{i=1}^n\|A^{(k)T}(X_i-\bar X)B^{(k)}\|_F^2$ and
$\frac{1}{n}\sum_{i=1}^n\|A^{(k+1)T}(X_i-\bar X)B^{(k+1)}\|_F^2$.
Output $\{\widehat A,\widehat B\}=\{A^{(k+1)},B^{(k+1)}\}$.
\end{itemize}

For any fixed $A^{(k)}$ or $B^{(k+1)}$, the optimization problems in
Steps~1 and 2 are the usual eigenvalue-problems of sizes $p$ and
$q$, respectively. Hence, $A^{(k+1)}$ and $B^{(k+1)}$ can be easily
obtained. Moreover, the algorithm ensures the quantity
$\frac{1}{n}\sum_{i=1}^n\|A^{(k)T}(X_i-\bar X)B^{(k)}\|_F^2$ to be
monotonically increasing as $k$ increases and, hence, the solution
must exist since $\frac{1}{n}\sum_{i=1}^n\|A^{(k)T}(X_i-\bar
X)B^{(k)}\|_F^2$ is bounded above by
$\frac{1}{n}\sum_{i=1}^n\|X_i-\bar X\|_F^2$. Because GLRAM can only
find a local maximum (depends on the chosen random initial
$A^{(0)}$), multiple random initials are suggested by Ye (2005) to
ensure the global maximum. In contrast to this suggestion, we
propose to use the leading $\tilde p$ eigenvectors of
$\frac{1}{n}\sum_{i=1}^n(X_i-\bar{X})(X_i-\bar{X})^T$ as an initial
of $A^{(0)}$.

We observe that hierarchical nesting structure may not exist for
MPCA. Precisely, if $(\ptilde',\qtilde')\le(\ptilde,\qtilde)$ with
the corresponding solution pairs $\{\widehat A',\widehat B'\}$ and
$\{\widehat A,\widehat B\}$, respectively, there is no guarantee
that $\rmspan(\widehat A')\subseteq\rmspan(\widehat A)$, nor
$\rmspan(\widehat B')\subseteq\rmspan(\widehat B)$. In the
population level, however, there certainly exist relationships
between the target subspaces and the MPCA subspaces prescribed by
the optimization problem~(\ref{max_prob}).

\begin{prop}\label{relationship}
Assume model (\ref{model}) and let $\{A,B\}$ be the solution pair to
the maximization problem (\ref{max_prob}) under dimensionality
$(\ptilde, \qtilde)$.
\begin{itemize}
\item[(a)] If $\ptilde\geq p_0$ and $\qtilde\geq q_0$, then
$\rmspan( A)\supseteq\rmspan(A_0)$ and $\rmspan(B)\supseteq\rmspan(
B_0)$.

\item[(b)] If $\ptilde < p_0$ and $\qtilde \geq q_0$, then
$\rmspan(A)\subsetneq\rmspan(A_0)$ and $\rmspan(B)\supseteq\rmspan(
B_0)$.

\item[(c)] If $\ptilde \geq p_0$ and $\qtilde < q_0$, then
$\rmspan(A)\supseteq\rmspan(A_0)$ and $\rmspan(B)\subsetneq\rmspan(
B_0)$.

\item[(d)] If $\ptilde < p_0$ and $\qtilde < q_0$, then
$\rmspan(A)\subsetneq\rmspan(A_0)$ and
$\rmspan(B)\subsetneq\rmspan( B_0)$.
\end{itemize}
\end{prop}

Even though there is no general hierarchical nesting structure for
MPCA subspaces, Proposition~\ref{relationship} ensures the existence
of a specific nesting structure, which the extracted MPCA subspace
is a proper subspace of the target subspace if the dimension is
under-specified, and contains the target subspace if the dimension
is over-specified. It also implies that MPCA indeed searches the
true target subspace $\rmspan(B_0\otimes A_0)$ when
$(\ptilde,\qtilde)=(p_0,q_0)$ is correctly specified. As a result,
these arguments provide a justification of using $\rmspan(\widehat
B\otimes \widehat A)$ in the sample level for subsequent statistical
analysis.

\subsection{Connection with ${\rm \textbf{(2D)}}^2$PCA and conventional PCA}

The ${\rm (2D)}^2$PCA is another method  to extract basis for tensor
objects. For a given dimensionality $(\ptilde,\qtilde)$, the
population ${\rm (2D)}^2$PCA components
$A^*=[a_1^*,\cdots,a_{\ptilde}^*]$ and
$B^*=[b_1^*,\cdots,b_{\qtilde}^*]$ are defined to be the leading
$\ptilde$ and $\qtilde$ eigenvectors of $E[(X-\mu)(X-\mu)^T]$ and
$E[(X-\mu)^T(X-\mu)]$ with the corresponding eigenvalues
$\{\lambda_i^*:1\leq i\leq \ptilde\}$ and $\{\xi_j^*:1\leq j\leq
\qtilde\}$. The sample analogues, denoted by $\widehat A^*$,
$\widehat B^*$, $\hat\lambda_i^*$, and $\hat\xi_j^*$ are similarly
defined to be the leading eigenvectors and eigenvalues of
$\frac{1}{n}\sum_{i=1}^n(X_i-\bar X)(X_i-\bar X)^T$ and
$\frac{1}{n}\sum_{i=1}^n(X_i-\bar X)^T(X_i-\bar X)$. The following
proposition states a connection between the ${\rm (2D)}^2$PCA and
MPCA in the population level.

\begin{prop}\label{thm.2d}
Assume model (\ref{model}) and  that $\{\lambda_i^*:1\leq i\leq
p_0\}$ and $\{\xi_j^*:1\leq j\leq q_0\}$ are simple roots.
\begin{itemize}
\item[(a)]
If $\,\qtilde\geq q_0$, then MPCA and (2D)$^2$PCA share the same
leading $(p_0\wedge \ptilde)$ eigenvectors, i.e., $a_i=a_i^*$ for
$i=1,\cdots,(p_0\wedge \ptilde)$. Moreover,
$\lambda_i^*-\lambda_i=(q-\qtilde)\sigma^2$ for $i=1,\dots,\ptilde$.
\item[(b)] If $\,\ptilde\geq p_0$, then MPCA and (2D)$^2$PCA share the same
leading $(q_0\wedge \qtilde)$ eigenvectors, i.e., $b_j=b_j^*$ for
$i=1,\cdots,(q_0\wedge \qtilde)$. Moreover,
$\xi_j^*-\xi_j=(p-\ptilde)\sigma^2$ for $j=1,\dots,\qtilde$.
\end{itemize}
\end{prop}

When the dimension $(\ptilde,\qtilde)$ is adequate,
Proposition~\ref{thm.2d} implies that ${\rm (2D)}^2$PCA and MPCA, in
the population level, actually target the same subspace
$\rmspan(B_0\otimes A_0)$ under model~(\ref{model}). However, there
is no guarantee that the extracted bases $\{\widehat A^*,\widehat
B^*\}$ of ${\rm (2D)}^2$PCA also maximize the sample version of
(\ref{max_prob}). Though, under the setting of
Proposition~\ref{thm.2d}, $E[(X-\mu)P_{B}(X-\mu)^T]$ and
$E[(X-\mu)(X-\mu)^T]$ have the same leading eigenvectors, we expect
an efficiency gain in using $E[(X-\mu)P_{B}(X-\mu)^T]$, since it is
less noise-contaminated than $E[(X-\mu)(X-\mu)^T]$.
A rigorous proof of efficiency gain is
provided in Section~3.


\begin{rmk}
From Proposition~\ref{thm.2d}, it is suggested to select
$(\ptilde,\qtilde)$ through ${\rm (2D)}^2$PCA, since the dimension
of $A_0$ and $B_0$ are not known before being estimated. A formal
statistical test is provided in Section 2.3.
\end{rmk}

There is also a connection between MPCA and  conventional PCA. Under
model (\ref{model}), without considering the random noise
$\varepsilon$, $\vec(X-\mu)$ belongs to $\rmspan(B_0 \otimes A_0)$,
which is the target subspace of MPCA. Observe also that
\begin{eqnarray}
\Sigma&=&(B_0\otimes A_0)T(B_0\otimes A_0)^T+\sigma^2 I_m \nonumber\\
&=&(B_0\otimes A_0)(T+\sigma^2 I_{m_0})(B_0\otimes A_0)^T +\sigma^2
Q_{B_0\otimes A_0}\label{Sigma_form},
\end{eqnarray}
where $Q_{B_0\otimes A_0}=Q_{B_0}\otimes P_{A_0}+P_{B_0}\otimes
Q_{A_0}+Q_{B_0}\otimes Q_{A_0}$ is the projection matrix onto the
complement of $\rmspan(B_0\otimes A_0)$. It should be noted that the
matrix $T$ is not necessarily a diagonal matrix. Hence, $(B_0\otimes
A_0)$ is not the same with the conventional PCA components in
general. If we further diagonalize $T=GDG^T$ with $D$ being a
diagonal matrix of size $m_0\times m_0$, we have the following
factorization:
\begin{eqnarray*}\label{cPCA}
\Sigma =[\Gamma,~\Gamma_\bot]
\left[\begin{array}{cc}D+\sigma^2 I_{m_0}&0\\
0&\sigma^2 I_{m-m_0}\end{array}\right] [\Gamma,~\Gamma_\bot]^T,
\end{eqnarray*}
where $\Gamma=(B_0\otimes A_0)G$ and $\Gamma_\bot$ is an
orthonormal basis for the orthogonal complement of
$\rmspan(\Gamma)$. Consequently, the conventional PCA uses
$\Gamma=(B_0\otimes A_0)G$ as coordinate system for a compressed
representation for $\vec(X-\mu)$, while the MPCA uses $(B_0\otimes
A_0)$. Notice that $\rmspan(B_0\otimes A_0)=\rmspan(\Gamma)$
provided that $T$ in (\ref{Sigma_form}) is of full rank. In
summary, MPCA and conventional PCA use the same subspace for
compressed data representation. However, MPCA requires less
parameters (see the following remark) to specify the
low-dimensional subspace than the conventional approach.

\begin{rmk}\label{correct_dim}
The number of free parameters required for  MPCA is $p_0 p-\frac12
p_0(p_0+1) +q_0 q-\frac12 q_0(q_0+1)$, which is relatively small in
contrast  with the number of free parameters required for
conventional PCA: $p_0q_0 pq-\frac12 p_0q_0 (p_0q_0+1)$. It is the
adoption of $(B_0\otimes A_0)$ for the sake of parsimony, which is
one of the purposes of using MPCA. The following table gives the
numbers of parameters needed to specify an orthonormal basis for a
subspace of dimensionality $p_0\times q_0$ within a space of
dimensionality $p\times q=100$. We fix $(p,q,p_0)=(10,10,5)$ and let
$q_0$ vary.
\begin{table}[h]
\begin{center}
\begin{tabular}{cccccc}
\hline \hline
$q_0$   &   1   &   2   &   3   &   4   &   5   \\
\hline
MPCA    &   44  &   52  &   59  &   65  &   70  \\
PCA &   485 &   945 &   1380    &   1790    &   2175    \\
\hline\\
\end{tabular}
\caption{\label{t0} \sf Numbers of required free parameters at
$(p,q,p_0)=(10,10,5)$.}
\end{center}
\end{table}
\end{rmk}

We remind the reader that there is no obvious ordering
relationship between the MPCA components and conventional PCA
components. This can be seen in a simple example when $T={\rm
Cov}(\vec(U))={\rm diag}(\vec(C))$, where $C$ is a matrix with
$C_{ij}={\rm Var}(U_{ij})$. For the case of uncorrelated
$U_{ij}$'s, $T$ is diagonal, and hence, the conventional PCA and
the MPCA share the same eigenvectors. The leading $p_0q_0$
eigenvalues of the conventional PCA are $\{C_{ij}+\sigma^2:1\leq
i\leq p_0,1\leq j\leq q_0\}$, which have a natural ordering
depending on the values of $C_{ij}$'s. On the other hand, the
leading eigenvalues of MPCA at $(p_0,q_0)$ are derived to be
$C_{i\centerdot}=\sum_{j=1}^{q_0}C_{ij}$, $i=1,\cdots,p_0$, and
$C_{\centerdot j}=\sum_{i=1}^{p_0}C_{ij}$, $j=1,\cdots,q_0$, where
the ordering depends on the column sums and row sums of
$C_{ij}$'s. Therefore, even if we pick $a_i$ and $b_j$ from
leading eigenvectors of $A$ and $B$, there is no guarantee that,
when paired together, $b_j\otimes a_i$ is on the top list of
leading eigenvectors of the conventional PCA.

\subsection{Selection of dimensionality}

This section is devoted to the selection of the dimensionality
$(\ptilde,\qtilde)$. Similar to the conventional PCA, we propose
that the dimension is determined by the explained variance, as a
popular method in conventional PCA. First we define the cumulative
variance, which is a measure of the total variance of the tensor
objects projected onto MPCA subspace.

\begin{defn}
Let $\{A,B\}$ be a solution pair to the problem~(\ref{max_prob}).
We call the quantity $\Phi(\ptilde,\qtilde)= E\|A^T(X-\mu)B\|_F^2$
the cumulative variance for $X$ at rank-$(\ptilde,\qtilde)$, and
the quantity
\begin{equation*}\label{ratio_eigenval}
\rho(\ptilde,\qtilde)=\frac{\Phi(\ptilde,\qtilde)}{\Phi(p,q)}
\end{equation*}
the explained percentage of total variance of $X$ at
rank-$(\ptilde,\qtilde)$. Note that $\Phi(p,q)=E\|X-\mu\|_F^2$.
The corresponding sample analogues are defined to be
$\widehat\Phi(\ptilde,\qtilde)= \frac{1}{n}\sum_{i=1}^n\|\widehat
A^T(X_i-\bar X)\widehat B\|_F^2$, $\widehat\Phi(p,q)=
\frac{1}{n}\sum_{i=1}^n\|X_i-\bar X\|_F^2$ and
\begin{equation*}\label{s_ratio_eigenval}
\hat\rho(\ptilde,\qtilde)=\frac{\widehat\Phi(\ptilde,\qtilde)}{\widehat\Phi(p,q)}.
\end{equation*}
\end{defn}

\begin{rmk}
Note that $\Phi(\ptilde,\qtilde)>\Phi(\ptilde',\qtilde')$ does not
necessarily imply $\rmspan(A_{\ptilde}\otimes B_{\qtilde})\supset
\rmspan(A_{\ptilde'}\otimes B_{\qtilde'})$. Similar phenomenon can
be observed on the cumulative distribution function. For instance,
in a 2-dimensional c.d.f $F$, the phenomenon
``$F(x_1,x_2)>F(x_1',x_2')$'' does not imply $\{(u,v):u\le
x_1,v\le x_2\}\supset\{(u,v):u\le x_1',v\le x_2'\}$.
\end{rmk}

From the description below Definition~\ref{def}, we have
$\Phi(\ptilde,\qtilde)
=\sum_{i=1}^{\ptilde}\lambda_i=\sum_{j=1}^{\qtilde}\xi_j$ and
$\widehat\Phi(\ptilde,\qtilde)=\sum_{i=1}^{\ptilde}\hat\lambda_i
=\sum_{j=1}^{\qtilde}\hat\xi_j$. Note that $\lambda_i$ and $\xi_j$,
as well as $\hat\lambda_i$ and $\hat\xi_j$, depend on the specified
dimensionality $(\ptilde,\qtilde)$. Also note that
$\Phi(\ptilde,\qtilde)\leq \Phi(p,q)$ always holds. Thus,
$\rho(\ptilde,\qtilde)\leq 1$ and is used as a measure of adequacy
for MPCA at dimensionality $(\ptilde,\qtilde)$. Specifically, for a
given $\rho_0\in (0,1)$, consider the hypothesis test:
\begin{equation}\label{test}
H_0: \rho(\ptilde,\qtilde)\le\rho_0 ~~~v.s.~~~H_1:
\rho(\ptilde,\qtilde)>\rho_0.
\end{equation}
A rejection of $H_0$ then indicates the chosen dimensionality
$(\ptilde,\qtilde)$ satisfies the condition that
$\rho(\ptilde,\qtilde)$ reaches the required level of explained
variance at a certain confidence. To perform the test, a reference
distribution for the sample analogue $\hat\rho(\ptilde,\qtilde)$
is required. We derive the asymptotic distribution of
$\sqrt{n}(\hat\rho(\ptilde,\qtilde)-\rho(\ptilde,\qtilde))$ in
Section 3, which can be used to construct the rejection region of
the test.


\section{Asymptotic properties for MPCA}

In this section, we investigate the asymptotic behavior of MPCA.
Without loss of generality, we assume $\mu=0$ to simplify the
notations in the rest of discussion. It then implies
$\Sigma=E[\vec(X)\vec(X)^T]$ and the population kernel matrices of
MPCA at dimensionality $(\ptilde,\qtilde)$ can be simplified to be
$E[XP_BX^T]$ and $E[X^TP_AX]$. Note also that the population kernel
matrices of (2D)$^2$PCA reduce to $E[XX^T]$ and $E[X^TX]$ in this
situation.

Let $S_n$ be the sample covariance matrix of
$\{\vec(X_i)\}_{i=1}^n$, where $X_i$'s are iid observations with
finite second moments following model (\ref{model}). By the central
limit theorem, we have
\begin{eqnarray}
\sqrt{n}(S_n-\Sigma)\stackrel{d}{\rightarrow}N,
\end{eqnarray}
where $\vec(N)$ is an $m^2$-variate normal with zero-mean and
covariance matrix $\Sigma_N=\cov
\left(\vec(X)\otimes\vec(X)\right)$. If $\vec(X)$ is further assumed
to be normally distributed, then $S_n$ follows a Wishart
distribution and $\Sigma_N$ is derived to be (Anderson, 1963)
\begin{eqnarray}
\Sigma_N=(I_{m^2}+K_{m,m})(\Sigma\otimes\Sigma),\label{normality}
\end{eqnarray}
where $K_{\ell,k} =\sum_{i=1}^\ell \sum_{j=1}^k H_{ij}\otimes
H_{ij}^T$ is the commutation matrix, and $H_{ij}$ is an
$\ell\times k$ matrix with one in the $(i,j)^{\rm th}$ entry and
zeros elsewhere. Some important properties involving commutation
matrix are listed here (Magnus and Neudecker, 1979). Let $M_1\in
\mathfrak{R}^{a_1\times b_1}$ and $M_2\in \mathfrak{R}^{a_2\times
b_2}$ be two arbitrary matrices. Then,
$K_{a_1,b_1}=K_{b_1,a_1}^T$, $K_{a_1,b_1}K_{b_1,a_1}=I_{a_1b_1}$,
$K_{a_1,b_1}=I_{a_1}$ if $b_1=1$,
$\vec(M_1^T)=K_{a_1,b_1}\vec(M_1)$, and $(M_2\otimes
M_1)=K_{a_2,a_1}(M_1\otimes M_2)K_{b_1,b_2}$. These properties
will be repeatedly used in the discussion of asymptotic theory
without further reference. We note that, unless explicitly
specified, the asymptotic properties derived in this section does
not rely on the normality of $\vec(X)$.

\subsection{Asymptotic distributions for principal components, projections,
cumulative variance and explained variance in MPCA}

We first state the weak convergence of the cumulative variances and
the tensor principal components of MPCA. The limiting distributions
for projections and explained variance are direct applications of
delta method.

\begin{thm}\label{thm.tpc}
Assume model (\ref{model}) and, for any fixed $(\ptilde,\qtilde)$
with $\ptilde \leq p_0$ and $\qtilde \leq q_0$, the leading
$\ptilde$ eigenvalues $\lambda_i(\Sigma,\ptilde,\qtilde)$'s and the
leading $\qtilde$ eigenvalues $\xi_j(\Sigma,\ptilde,\qtilde)$'s of
MPCA are simple roots.
\begin{itemize}
\item[(a)]
For $\ptilde \leq p_0$ and $\qtilde \leq q_0$,\footnote{For
$\ptilde>p_0$ (or $\qtilde>q_0$, resp.) $E[XX^T]$ (or $E[X^TX]$,
resp.) has multiple roots from the $(p_0+1)^{\rm th}$ (or
$(q_0+1)^{\rm th}$, resp.) eigenvalue and beyond.} we have the
limiting distribution
\begin{eqnarray}
\sqrt{n}\left(\left[\begin{array}{c}\widehat\Phi(\ptilde,\qtilde)\\
\widehat\Phi(p,q)
\end{array}\right]-\left[\begin{array}{c}\Phi(\ptilde,\qtilde)\\
\Phi(p,q)
\end{array}\right]\right)\stackrel{d}{\rightarrow}\left[\begin{array}{c}D_{\Phi(\ptilde,\qtilde)}\\
\vec(I_m)^T
\end{array}\right]\vec(N),
\end{eqnarray}
where $D_{\Phi(\ptilde,\qtilde)}=
\frac{\partial\Phi(\ptilde,\qtilde)} {\partial\vec(\Sigma)}$ and
its explicit expression is given in Lemma~\ref{lem.diff}.

\item[(b)]
For $\ptilde \leq p_0$ and $\qtilde \leq q_0$,\footnote{For either
$\ptilde>p_0$, or $\qtilde>q_0$, the $(p_0+1)^{\rm th}$, or
$(q_0+1)^{\rm th}$, tensor principal components are not uniquely
determined due to multiple characteristic roots.} we have the
limiting distribution
\begin{eqnarray}
\sqrt{n}\left(\left[\begin{array}{c}\vec(\widehat A)\\
\vec(\widehat B)
\end{array}\right]-\left[\begin{array}{c}\vec(A)\\
\vec(B)
\end{array}\right]\right)\stackrel{d}{\rightarrow}D_{H_{\ptilde,\qtilde}}\vec(N),
\end{eqnarray}
where
\begin{equation*}\label{DH}
D_{H_{\ptilde,\qtilde}}=\left[(\frac{\partial a_1}{\partial
\vec(\Sigma)})^T,\cdots,(\frac{\partial a_{\ptilde}}{\partial
\vec(\Sigma)})^T,(\frac{\partial b_1}{\partial
\vec(\Sigma)})^T,\cdots,(\frac{\partial b_{\qtilde}}{\partial
\vec(\Sigma)})^T\right]^T.
\end{equation*}
When $(\ptilde,\qtilde)=(p_0,q_0)$, $D_{H_{p_0,q_0}}$ has an
explicit expression, which is given in Lemma~\ref{lem.diff}.
\end{itemize}
\end{thm}

\begin{lem}\label{lem.diff}
Assume the model~(\ref{model}).
\begin{itemize}
\item[(a)]
For $\ptilde \leq p_0$ and $\qtilde \leq q_0$, we have
\begin{eqnarray}
D_{\Phi(\ptilde,\qtilde)}=\vec(P_{B\otimes A})^T.
\end{eqnarray}
\item[(b)]  When $(\ptilde,\qtilde)=(p_0,q_0)$, for $i=1,\cdots,p_0$ and
$j=1,\cdots,q_0$, we have
\begin{eqnarray}
\frac{\partial a_i}{\partial\vec(\Sigma)}&=&\left\{a_i\otimes
\vec(P_{B_0})\otimes (\lambda_iI_p-E[XP_{B_0}X^T])^{+}\right\}^T(K_{p,q}\otimes I_{pq})\\
\frac{\partial b_j}{\partial\vec(\Sigma)}&=&\left\{b_j\otimes
\vec(P_{A_0})\otimes
(\xi_jI_q-E[X^TP_{A_0}X])^{+}\right\}^T(I_{pq}\otimes K_{p,q}),
\end{eqnarray}
where, for a given matrix $M$, $M^+$ denotes its Moore-Penrose
generalized inverse.
\end{itemize}
\end{lem}

It can be seen from Lemma~\ref{lem.diff} that, when
$(\ptilde,\qtilde)=(p_0,q_0)$, the asymptotic distribution of
$\widehat A$ depends on $B$ only through $\rmspan(B)=\rmspan(B_0)$,
and the asymptotic distribution of $\widehat B$ depends on $A$ only
through $\rmspan(A)=\rmspan(A_0)$. We are now on the position to
obtain the asymptotic normality of the projection matrix onto MPCA
subspace $P_{\widehat B\otimes\widehat A}$ and the explained
variance $\hat\rho(\ptilde,\qtilde)$ in the following corollaries.

\begin{cor}\label{asymp.proj}
Under the same assumptions of Theorem~\ref{thm.tpc}. For $\ptilde
\leq p_0$ and $\qtilde \leq q_0$, we have the limiting distribution
of the projection matrix onto MPCA subspace
\begin{eqnarray}
\sqrt{n}~\vec(P_{\widehat B\otimes \widehat A}-P_{B\otimes
A})\stackrel{d}{\rightarrow}D_{P_{B\otimes A}}\vec(N),
\end{eqnarray}
where $D_{P_{B\otimes A}}=\frac{\partial \vec(P_{B\otimes
A})}{\partial \vec(\Sigma)}$. When $(\ptilde,\qtilde)=(p_0,q_0)$,
$D_{P_{B_0\otimes A_0}}$ has the explicit expression
\begin{eqnarray}
&&(I_{m^2}+K_{m,m})\times \nonumber\\
&&\Bigg\{\sum_{i=1}^{p_0} (K_{q,p}\otimes
I_{pq})\left(P_{a_i}\otimes[\vec(P_{B_0})\vec(P_{B_0})^T]\otimes
\{\lambda_iI_p-E[XP_{B_0}X^T]\}^{+}\right)(K_{p,q}\otimes I_{pq})\nonumber\\
&&+\sum_{j=1}^{q_0}(I_{pq}\otimes K_{q,p})
\left(P_{b_j}\otimes[\vec(P_{A_0})\vec(P_{A_0})^T]\otimes
\{\xi_jI_q-E[X^TP_{A_0}X]\}^{+}\right)(I_{pq}\otimes
K_{p,q})\Bigg\}.\nonumber
\end{eqnarray}
\end{cor}

\begin{cor}\label{asymp.ratio}
Under the same assumptions of Theorem~\ref{thm.tpc}. For $\ptilde
\leq p_0$ and $\qtilde \leq q_0$, we have the limiting distribution
of the explained variance
\begin{eqnarray}
\sqrt{n}(\hat\rho(\ptilde,\qtilde)-\rho(\ptilde,\qtilde))\stackrel{d}{\rightarrow}N(0,\sigma^2_{\rho(\ptilde,\qtilde)}),
\end{eqnarray}
where $\sigma^2_{\rho(\ptilde,\qtilde)}$ is defined to be
\begin{equation}
\left(\frac{1}{{\Phi(p,q)}}D_{\Phi(\ptilde,\qtilde)}-\frac{{\Phi(\ptilde,\qtilde)}}{{\Phi(p,q)}^2}\vec(I_m)^T\right)
\Sigma_N
\left(\frac{1}{{\Phi(p,q)}}D_{\Phi(\ptilde,\qtilde)}-\frac{{\Phi(\ptilde,\qtilde)}}{{\Phi(p,q)}^2}\vec(I_m)^T\right)^T.
\end{equation}
\end{cor}

Corollary~\ref{asymp.ratio} is the cornerstone of our asymptotic
test for hypothesis~(\ref{test}). Before practical implementation
of the test, however, we need a consistent estimator of
$\sigma^2_{\rho(\ptilde,\qtilde)}$. Note that the asymptotic
covariance $\Sigma_N$ can be empirically estimated by
\[\widehat\Sigma_{N,1}=\frac{1}{n}\sum_{i=1}^n\vec\left(\vec(X_i-\bar
X)\vec(X_i-\bar X)^T-S_n\right) \vec\left(\vec(X_i-\bar
X)\vec(X_i-\bar X)^T-S_n\right)^T.\] Moreover, if $\vec(X)$ is
normally distributed, we can also estimate $\Sigma_N$ by
\[\widehat\Sigma_{N,2}=(I_{m^2}+K_{m,m})(S_n\otimes S_n)\] based
on~(\ref{normality}). Consequently, the asymptotic variance
$\sigma_{\rho(\ptilde,\qtilde)}^2$ is estimated by
$$\widehat\sigma^2_{\rho(\ptilde,\qtilde)}=\left(\frac{1}{{\widehat\Phi(p,q)}}\widehat D_{\Phi(\ptilde,\qtilde)}-
\frac{{\widehat\Phi(\ptilde,\qtilde)}}{{\widehat\Phi(p,q)}^2}\vec(I_m)^T\right)
\widehat\Sigma_{N,i} \left(\frac{1}{{\widehat\Phi(p,q)}}\widehat
D_{\Phi(\ptilde,\qtilde)}-
\frac{{\widehat\Phi(\ptilde,\qtilde)}}{{\widehat\Phi(p,q)}^2}\vec(I_m)^T\right)^T$$
for $i=1,2$ (depends on the normality of $\vec(X)$ or not), where
$\widehat D_{\Phi(\ptilde,\qtilde)}=\vec(P_{\widehat B
\otimes\widehat A})^T$. The consistency of
$\widehat\sigma_{\rho(\ptilde,\qtilde)}^2$ is a direct consequence
by standard arguments. These facts enable us to construct an
approximate level $\alpha$ test to determine the dimensionality
$(\ptilde,\qtilde)$.

\begin{thm}\label{asymp.test}
Assume the conditions of Theorem~\ref{thm.tpc} and
$(\ptilde,\qtilde)\leq(p_0,q_0)$. For the hypothesis (\ref{test}),
an approximated level $\alpha$ test is to reject $H_0$ if
\begin{equation}
\hat\rho(\ptilde,\qtilde)
> \rho_0+\frac{\hat\sigma_{\rho(\ptilde,\qtilde)}}{\sqrt{n}}z_\alpha,
\end{equation}
where $z_\alpha$ is the upper $\alpha$ quantile of the standard
normal.
\end{thm}

\subsection{Asymptotic efficiency}

MPCA and ${\rm (2D)}^2$PCA actually target the same basis when
$(\ptilde,\qtilde)=(p_0,q_0)$. Intuitively, we are in favor of MPCA
since its kernel matrices are less noise-contaminated than the ones
of ${\rm (2D)}^2$PCA as mentioned previously. The following theorem
proves that MPCA is indeed asymptotically more efficient than ${\rm
(2D)}^2$PCA, wherein $\acov$ denotes the asymptotic covariance.

\begin{thm}\label{thm.efficiency}
Assume the conditions of Theorem~\ref{thm.tpc} and the normality of
$\vec(X)$. Let $(\ptilde,\qtilde)=(p_0,q_0)$ and let $(\widehat
A^*,\widehat B^*)$ be the ${\rm (2D)}^2$PCA components under
$(p_0,q_0)$. Then,
\begin{eqnarray}\label{eq34}
{\rm aCov}(\vec(P_{\widehat B^*\otimes \widehat A^*}))-{\rm
aCov}(\vec(P_{\widehat B\otimes \widehat A}))\geq 0,
\end{eqnarray}
where the equality holds if and only if $(p_0,q_0)=(p,q)$.
\end{thm}

Theorem~\ref{thm.efficiency} states that under
model~(\ref{model}), MPCA is at most as disperse as (2D)$^2$PCA in
estimating the dimension reduction subspace $\rmspan(B_0\otimes
A_0)$. The only case that we will gain nothing from MPCA over the
(2D)$^2$PCA is when $(p_0,q_0)=(p,q)$.  Note that the condition
$(p_0,q_0)=(p,q)$ implies that there is no room of dimension
reduction at all and is probably of no interest in real
applications. Consequently, Theorem~\ref{thm.efficiency} provides
a justification of using MPCA.


\section{Experimental study: the {\sf Olivetti Faces} data set}

We test and compare the performance of MPCA and conventional PCA on
{\sf Olivetti Faces} data set, which is available at {\tt
http://www.cs.nyu.edu/$\sim$roweis/data.html}. This data set
consists of 400 gray scale (8 bits) face images of  $\,64\times 64$
pixels. There exist different facial expressions and/or views for
each individual in this data set. A simulation experiment is
designed as follows. 400 face images are randomly partitioned into a
training set with size 100 and a test set with size 300. This
\mbox{100-300} partition, where the training set is smaller than
test set, is to reflect a scenario of using a small portion of data
to train a basis set for the representation of the rest data in data
archive.

Both MPCA and conventional PCA are applied on the 100 training
images to produce image basis which is used to reconstruct the rest
300 test images. The average of the 100 training images, named mean
face, has been subtracted from all the 400 images for PCA training
and for test image reconstructions as well. The mean face is finally
added to the reconstructions at the last stage to show the resulting
images. 500 replicates of training-test partitions are performed to
compare the mean test error, which is defined as the average of the
Frobenius norm between the original images and the reconstructed
images on test data set. The result is in Table~\ref{t1}. The mean
test error for conventional PCA is more than seven times of that for
MPCA; and the standard deviation is more than 12 times.

\begin{table}[h]
\begin{center}
\smallskip
\begin{tabular}{|c|cc|}\hline
 Frobenius-Norm & MPCA &  Conventional PCA \\ \hline
 Mean $(\times 10^5)$ &   1.1346    &  8.6455\\      \hline
 SD   $(\times 10^2)$&    9.6398  & 120.39     \\   \hline
\end{tabular}
\caption{\label{t1}\sf Test error comparison for MPCA and
conventional PCA on the test images of {\sf Olivetti Faces} data set.
The error is defined as the Frobenius norm of two image matrices:
original test image and its reconstruction using $28\times 28$
principal components. }
\end{center}
\end{table}

\begin{figure}[h]
\vspace{-0.5cm} \hspace{-2.8cm}
\includegraphics[height=13cm]{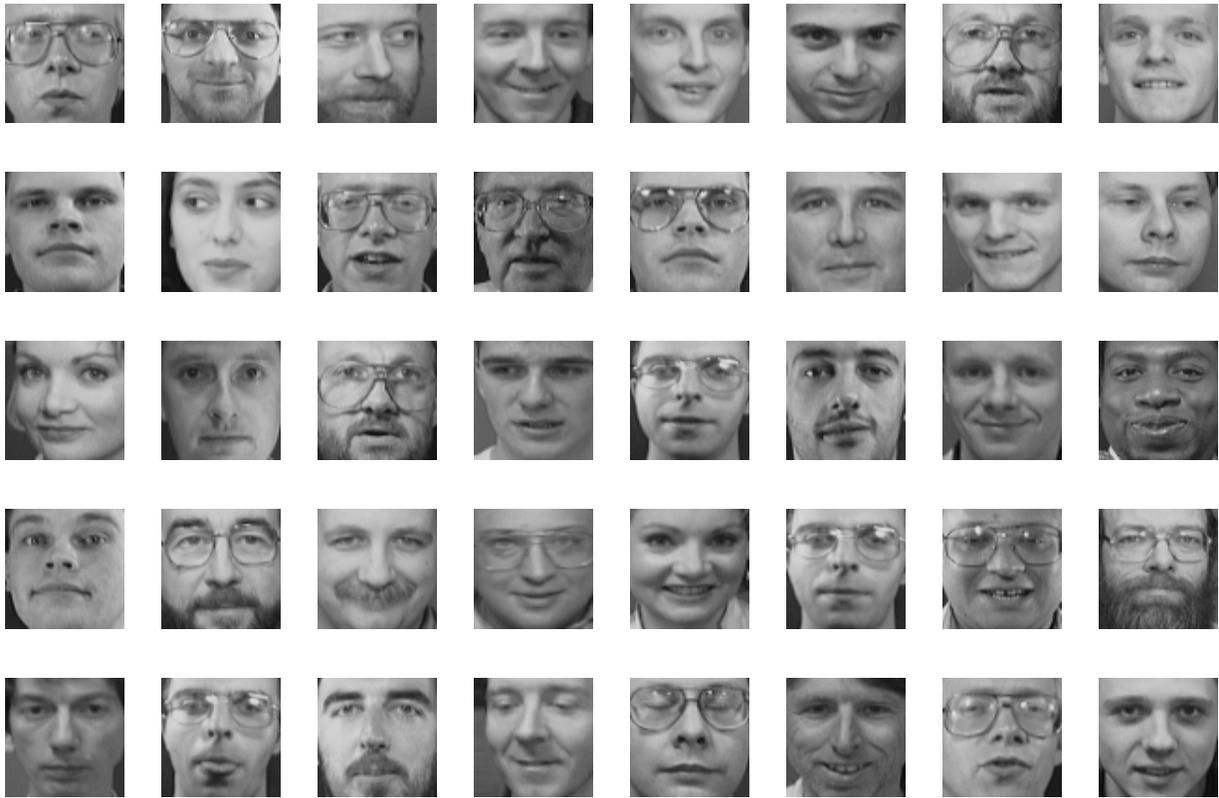}
\vspace{-2.5cm} \caption{\label{test.face} \sf 40 test faces
randomly drawn from the test set.}
\end{figure}

\begin{figure}[h]
\vspace{-0.5cm} \hspace{-2.8cm}
\includegraphics[height=13cm]{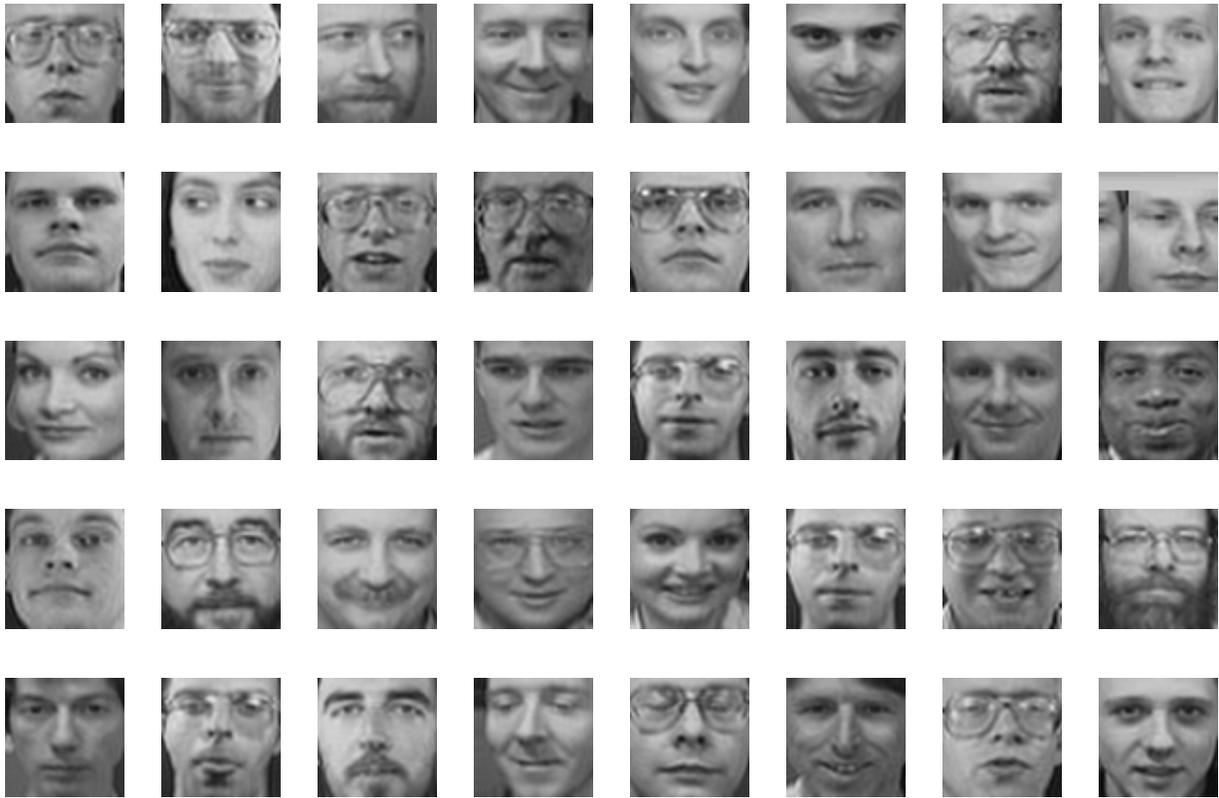}
\vspace{-2.5cm} \caption{\label{recon.tensor}\sf Reconstructed faces
using $28\times 28$ trained MPCA components.}
\end{figure}

\begin{figure}[h]
\vspace{-0.5cm} \hspace{-2.8cm}
\includegraphics[height=13cm]{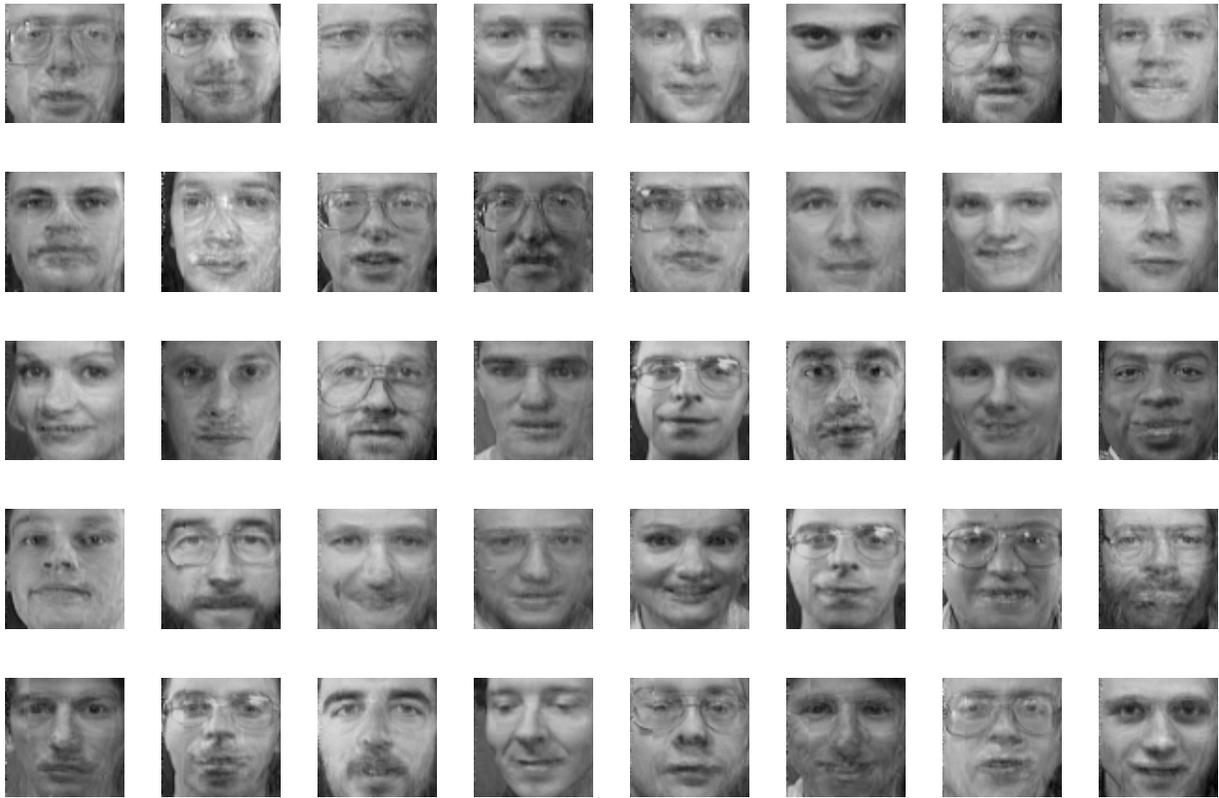}
\vspace{-2.5cm} \caption{\label{recon.conven} \sf Reconstructed
faces using $784$ trained conventional PCA components.}
\end{figure}

\begin{figure}[h]
\vspace{0cm} \hspace{0cm}
\includegraphics[height=13cm]{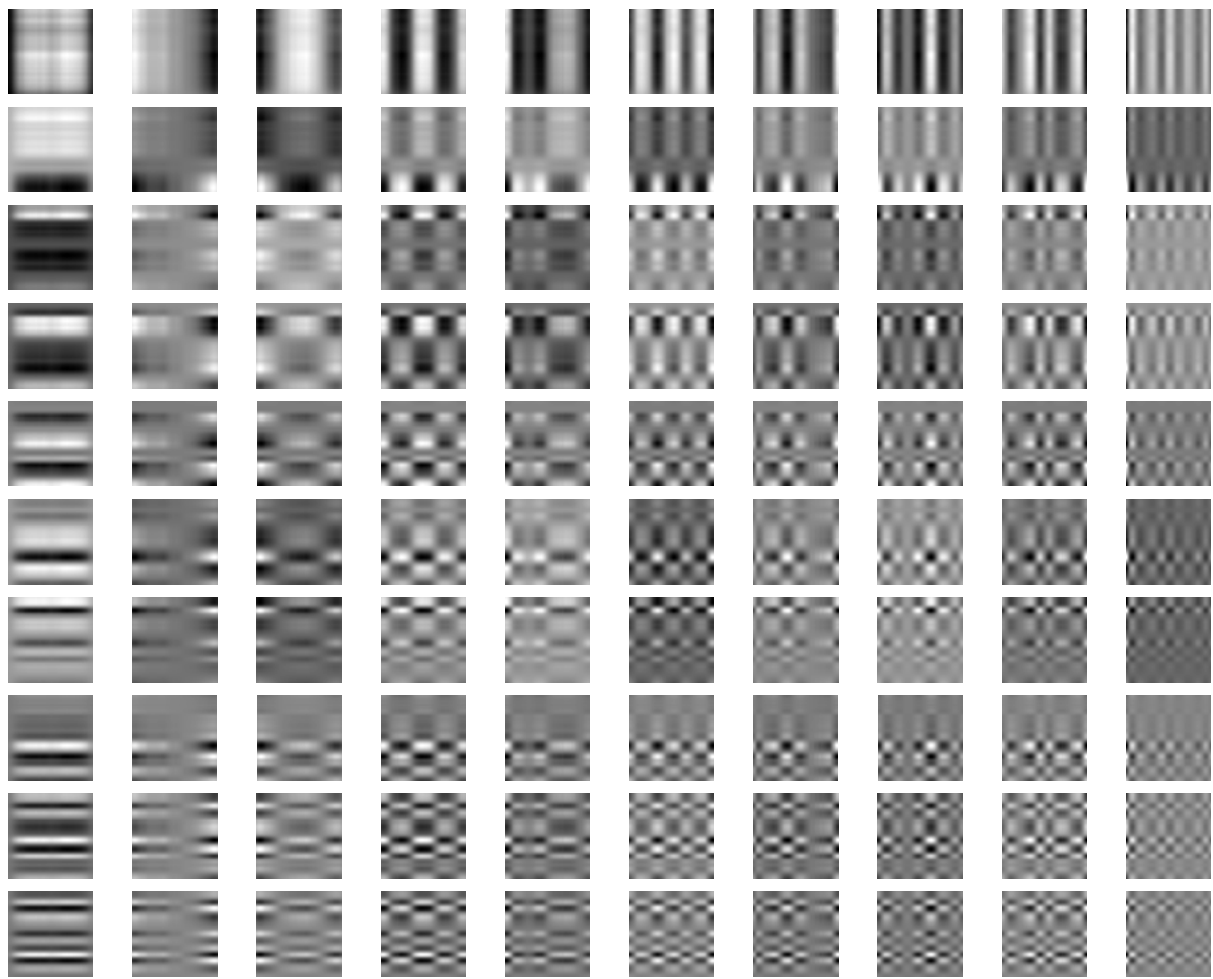}
\vspace{-1cm} \caption{\label{basis.tensor} \sf Leading 100 MPCA
basis images.}
\end{figure}

\begin{figure}[h]
\includegraphics[width=1\textwidth]{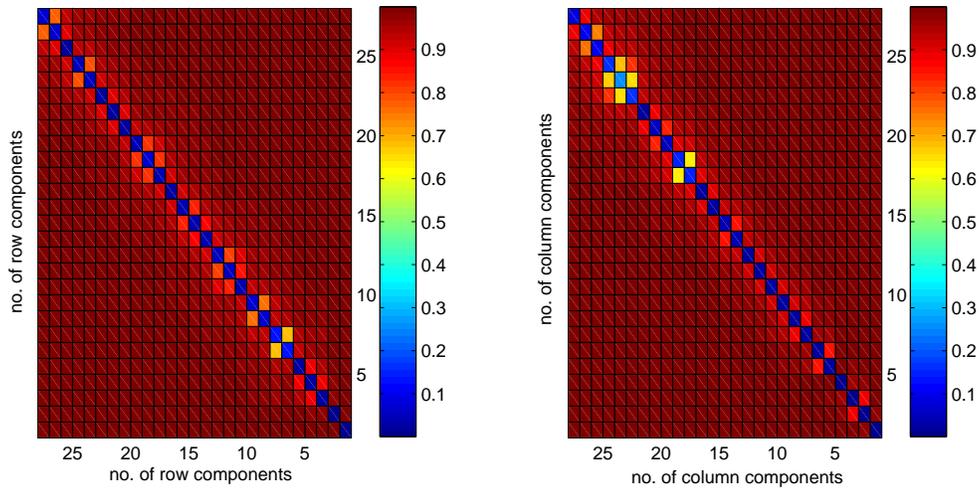}
\caption{\label{v-plot}\sf The variability pattern plots for MPCA
for $1\le\ptilde\le 28$ and $1\le\qtilde\le 28$.}
\end{figure}

In Figures~\ref{test.face}-\ref{recon.conven}, 40 test images are
randomly chosen from the test set to show the visual performance of
image reconstructions by these two PCA schemes. In MPCA, 28 row
eigenvectors and 28 column eigenvectors, both with size 64, are used
to generate 784 basis images, of which the 100 leading ones are
shown in Figure~\ref{basis.tensor}. We remind the reader that the
selection $(\ptilde,\qtilde)=(28,28)$ produces an
$\hat\rho(\ptilde,\qtilde)$ value $0.968$. Based on
Theorem~\ref{asymp.test}, a one-sided $95\%$ confidence interval for
$\hat\rho(28,28)$ is given by $[0.967,1]$. We also show the
variability pattern plots (Tu and Huang, 2011) in Figure
\ref{v-plot}. These plots present the average variations (absolute
values) of the eigenvectors for the bootstrap re-sampled data, from
those eigenvectors for the original data. The horizontal and
vertical indices refer to the eigenvector indices for re-sampled and
original data. The indices of eigenvectors are sorted by
eigenvalues. The variations are presented by colors from dark blue
for perfect matched, to dark red for extremely  deviated. Usually,
eigenvectors with distinct eigenvalues show deep blue on the
diagonal and deep red on the off-diagonal. Eigenvectors with the
same multiple root eigenvalue tend to be visualized by a cubic
pattern on their correspondence indices. It can be seen that our
choice of $(\ptilde,\qtilde)=(28,28)$ does not produce multiple
roots, since the bootstrapped variability of the solutions at this
selection is quite small.

\begin{figure}[h]
\vspace{0cm} \hspace{0cm}
\includegraphics[height=13cm]{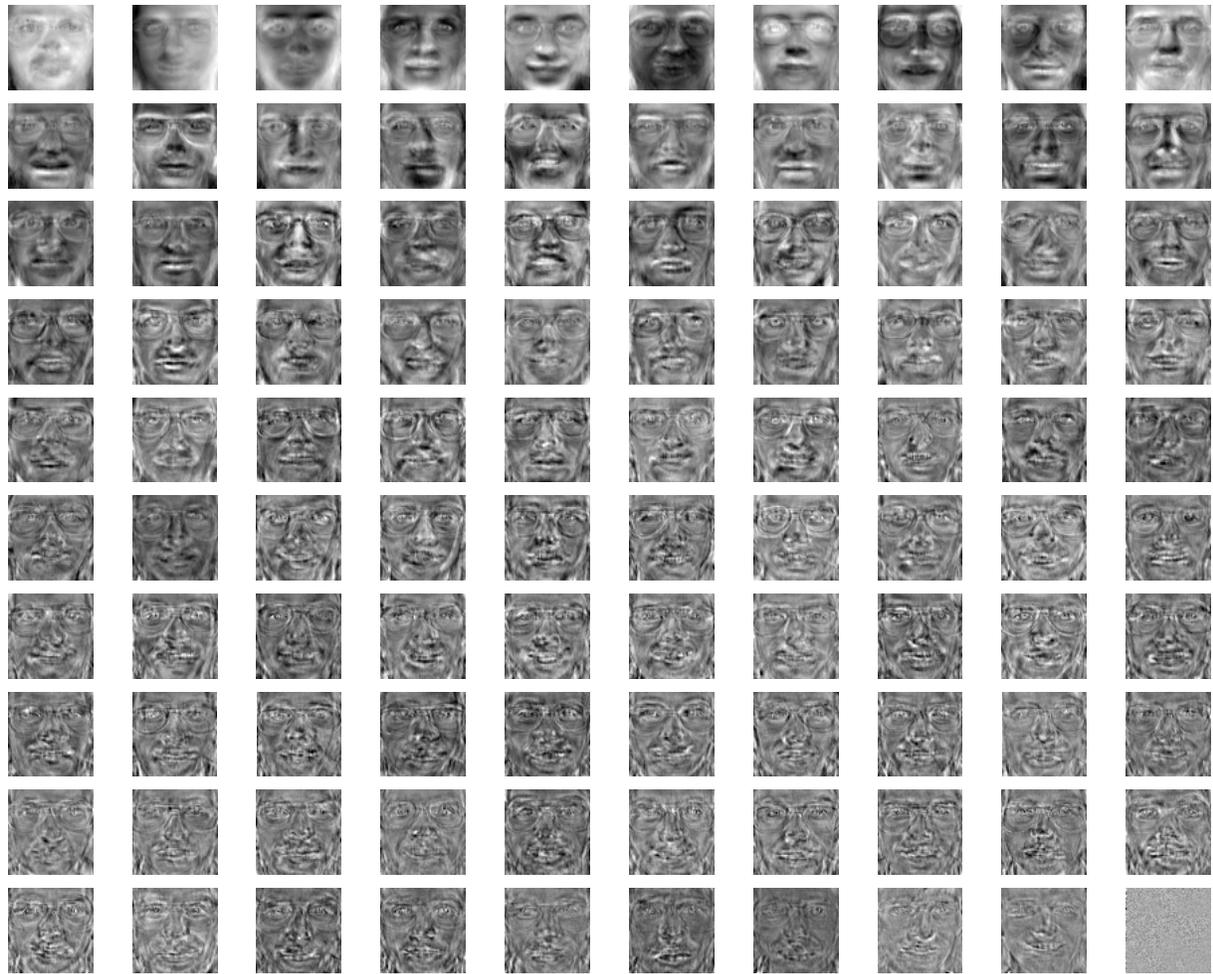}
\vspace{-1cm} \caption{\label{basis.conven}\sf Leading 100
conventional PCA basis images.}
\end{figure}

In conventional PCA, $784~(=28\times 28)$ eigenvectors (basis
images) with size 4096 are used, of which the 100 leading ones are
shown in Figure~\ref{basis.conven}. Because of using 100 training
images with average subtraction, there are at most 99 meaningful
eigenvectors in the conventional PCA. The rest are randomly
orthogonal eigenvectors with zero eigenvalue from the remaining
subspace. In Figure~\ref{basis.conven}, from top to bottom, we can
see the images with clear facial shape to vague ones and a random
image on the $100^{\rm th}$ one. On the other hand, MPCA tends to
distribute the image characteristics to more basis elements which
may allow for more local modification on the images.

\begin{figure}[h]
\begin{center}
\vspace{0cm}
\includegraphics[height=4.5cm]{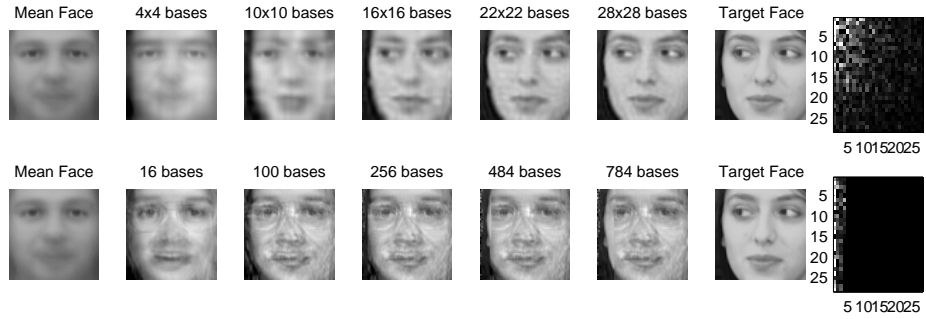}
\vspace{-.8cm}\caption{\label{increase_basis}\sf The reconstructed
images for the test face by adding more basis elements  are compared
for MPCA (top row) and conventional PCA (bottom row). Both mean face
and the target face are put in this figure as references. Plots in
the right-most column are projection coefficients (in absolute
value) onto PCA subspace.}
\end{center}
\end{figure}

In Figure~\ref{increase_basis}, one particular image among the 40
test images is chosen to demonstrate the performance of these two
methods. The top row shows the image reconstruction process for
MPCA when more basis elements are added in, and the bottom shows
for conventional PCA. The mean face is put in the first column and
the target image in the $7^{\rm th}$ column as references. The
right-most column shows the absolute values of projection scores
on the leading 784 basis elements. It is clear that the
conventional PCA concentrates on no more than 99 basis elements
while the MPCA spreads out to much more basis elements. For MPCA,
the image turns its view when $10\time 10$ basis elements are
used; the pupil turns to left when $16\times 16$ basis elements
are used; the double eyelid and nostrils show up when $22\times
22$ basis elements are used; the facial curves become clear when
$28\times 28$ basis elements are used. While we can observe the
reconstruction progress by adding more basis elements for MPCA, we
do not see much difference after 100 basis elements for
conventional PCA. It is clear that MPCA performs better than
conventional PCA in reconstructing the test images from
Table~\ref{t1} and these figures.

\section{Concluding discussions}

PCA is a  popular tool to reduce the dimensions for high dimensional
data analysis; MPCA could be likely to serve the similar function
for higher order tensor data sets. From this work, the statistical
properties of MPCA become clear through the theoretical framework
and the performance of MPCA is predictable through the asymptotic
results. Most importantly, based on these asymptotic results,
various hypothesis tests become feasible for subsequent analysis,
including pattern recognition or classification. Our work, though
technically theoretical, may construct a platform to expand the
application potentials of MPCA.

The advantages of MPCA over conventional PCA on tensor structure
data are evident in the {\sf Olivetti Faces} data example.
Therein, conventional PCA suffers seriously from the large $m$ and
small $n$ problem such that there can be at most $n-1$ meaningful
eigenvectors. This makes it unavoidable that all the data noises
are still carried by the chosen principal components. Furthermore,
too concentrated information in one component, which may not be
good for pattern recognition or classification prediction. On the
other hand, MPCA distributes the information to more components
which may allow local modification in the process of image
reconstruction, with even fewer free parameters. The key point for
the good performance of MPCA is the data tensor structure. For
practical purposes, the robustness of MPCA over model variety
should be further investigated.

\section*{References}

\begin{description}
\item
Anderson, T. W. (1963). Asymptotic theory for principal component
analysis. {\it Annals of Mathematical Statistics}, 34, 122-148.
\item
De Lathauwer, L., De Moor, B. and Vandewalle, J. (2000a). A
multilinear singular value decomposition. {\it SIAM J. Matrix Anal.
Appl.}, 21, 1253-1278.
\item
De Lathauwer, L., De Moor, B. and Vandewalle, J. (2000b). On the
best rank-1 and rank-$(R_1,R_2,\dots,R_N)$ approximation of
higher-order tensors. {\it SIAM J. Matrix Anal. Appl.}, 21,
1324-1342.
\item
Fine, J. (1987). On the validity of the perturbation method in asymptotic theory.
{\it Statistics}, 18, 401-414.
\item
Henderson, H. V. and Searle,  S. R. (1979). Vec and vech operators
for matrices, with some uses in Jacobians and multivariate
statistics. {\it Canadian J. Statistics}, 7, 65-81.
\item
Jolliffe, I.T. (2002). {\it Principal Component Analysis}.
Springer, New York.
\item
Kolda, T.G. and Bader, B.W. (2009). Tensor decompositions and
applications. {\it SIAM Review}, 51(3), 455-500.
\item
Li, B., Kim, M.K. and Altman, N. (2010). On dimension folding of
matrix- or array-valued statistical objects. {\it Annals of
Statistics}, 38, 1094-1121.
\item
Lu, H., Plataniotis, K. N. and Venetsanopoulos, A. N. (2008). MPCA:
Multilinear principal component analysis of tensor objects. {\it
IEEE Transactions on Neural Networks}, 19, 18-39.
\item
Magnus, J. R. and Neudecker, H. (1979). The commutation matrix: some
properties and applications. {\it Annals of Statistics}, 7, 381-394.
\item
Sibson, R. (1979). Studies in the robustness of multidimensional
scaling: perturbational analysis of classical scaling. {\it J.
Roy. Statist. Soc.}, 41, 217-229.
\item
Tu, I. P. and Huang, H. C. (2011). An estimation on a covariance
matrix when multiple roots exist. {\it manuscript}.
\item
Tyler, D. E. (1981). Asymptotic inference for eigenvectors. {\it
Annals of Statistics}, 9, 725-736.
\item
Yang, J., Zhang, D., Frangi, A.F. and Yang, J.Y. (2004).
Two-dimensional PCA: a new approach to appearance-based face
representation and recognition. {\it IEEE Transactions on Pattern
Analysis and Machine Intelligence}, 26, 131-137.
\item
Ye, J. (2005). Generalized low rank approximations of matrices. {\it
Machine Learning}, 61, 167-191.
\item
Zhang, D. and Zhou, Z. H. (2005). ${\rm (2D)}^2$PCA: Two-directional
two-dimensional PCA for efficient face representation and
recognition. {\it Neurocomputing}, 69, 224-231.
\end{description}

\section*{Appendix}
\renewcommand{\theequation}{A.\arabic{equation}} \setcounter{equation}{0}
\renewcommand{\thelem}{A.\arabic{lem}}

\begin{proof}[\textbf{Proof of Proposition~\ref{lrapprox}}]
In the maximization problem (\ref{max_prob}), the objective function
is continuous and the feasible region
$\O_{q,\qtilde}\otimes\O_{p,\ptilde}$ is compact. (Both continuity
and compactness are with respect to the topology induced by
Frobenius norm.) Thus, solution(s) exists.
\end{proof}

\begin{proof}[\textbf{Proof of Proposition~\ref{relationship}}]
(a) Let $[B,\,B_{\perp}]$ be a $q\times q$ orthonormal matrix. Since
$B_0\in\rmspan([B,\,B_{\perp}])$, there exists
$\eta_1\in\Re^{\qtilde\times q_0}$ and $\eta_2\in
\Re^{(q-\qtilde)\times q_0}$ such that
$B_0=B\eta_1+B_{\perp}\eta_2$. As $B_0^TB_0=I_{q_0}$, we have
$\eta_1^T\eta_1+\eta_2^T\eta_2=I_{q_0}$. Observe that
\begin{eqnarray}
&&E\|A^T(X-\mu)B\|_F^2\nonumber\\
&=&E\left\{{\tr}(A^T(A_0UB_0^T)BB^T(B_0U^TA_0^T)A)\right\}
+E\left\{\tr(A^T\varepsilon BB^T\varepsilon^TA)\right\}\nonumber\\
&=&E\left\{\tr(A^TA_0U\eta_1^T\eta_1 U^TA_0^TA)\right\}
+\ptilde\qtilde\sigma^2\nonumber\\
&=&E\left\{\tr(A^TA_0UU^TA_0^TA)\right\}-E\left\{\tr(A^TA_0U\eta_2^T\eta_2U^TA_0^TA)\right\}+\ptilde\qtilde\sigma^2\nonumber\\
&\leq&E\left\{\tr(A^TA_0UU^TA_0^TA)\right\}+\ptilde\qtilde\sigma^2,\label{eq.8}
\end{eqnarray}
where the equality in (\ref{eq.8}) holds if and only if $\eta_2=0$,
if and only if $\eta_1^T\eta_1=I_{q_0}$. Thus, if $\qtilde\geq q_0$,
such an $\eta_1$ (with rank $q_0$) exists to ensure the equality in
(\ref{eq.8}). This implies $B_0=B\eta_1$ and, hence,
$B_0\in\rmspan(B)$. Similarly, $A_0\in\rmspan(A)$ which establishes
(a).

To show (b), when $\qtilde\geq q_0$, from (a) we have
$B_0\in\rmspan(B)$ and
\begin{eqnarray}
\max_{A\in\O_{p,\ptilde},B\in\O_{q,\qtilde}}E\|A^T(X-\mu)B\|_F^2
&=&\max_{A\in\O_{p,\ptilde}} {\rm
trace}(A^TE[A_0UU^TA_0^T]A)+\ptilde\qtilde\sigma^2,
\end{eqnarray}
which is an eigenvalue-problem for the matrix $E[A_0UU^TA_0^T]$. By
diagonalizing $E[UU^T]=\Gamma_U\Lambda_U\Gamma_U^T$, then,
$E[A_0UU^TA_0^T]$ has $p_0$ non-zero eigenvalues $\Lambda_U$ with
the corresponding eigenvectors $A_0\Gamma_U$. When $\ptilde< p_0$,
the maximizer $A$ consists of the first $\ptilde$ columns of
$A_0\Gamma_U$ and, hence, $\rmspan( A)\subsetneq \rmspan(A_0)$.

(c) can be established in a similar way as (b).

To show (d), observe that
\begin{equation}\label{prop.3.4.d}
E\|A^T(X-\mu)B\|_F^2=E\left\{\tr(A^TA_0UB_0^TBB^TB_0U^TA_0^TA)\right\}
+\ptilde\qtilde\sigma^2.
\end{equation}
To maximize (\ref{prop.3.4.d}) over $A\in\O_{p,\ptilde},
B\in\O_{q,\qtilde}$ with $\ptilde < p_0$ and $\qtilde < q_0$, the
rank of $A^TA_0$ and $B^TB_0$ must be $\ptilde$ and $\qtilde$,
respectively, in order to attain the maximal value. This can happen
only if $\rmspan(A)\subsetneq\rmspan(A_0)$ and
$\rmspan(B)\subsetneq\rmspan(B_0)$.
\end{proof}

\begin{proof}[\textbf{Proof of Proposition~\ref{thm.2d}}]
We will only provide a proof for (a), and (b) can be obtained in a
similar way. If $\qtilde\geq q_0$, from
Proposition~\ref{relationship}~(a) we have
$\rmspan(B_0)\subseteq\rmspan(B)$, which further implies that
\begin{eqnarray}
E[(X-\mu)P_B(X-\mu)^T]&=&E[A_0UU^TA_0^T]
+E[\varepsilon P_B\varepsilon^T] \label{eig_A},\\
E[(X-\mu)(X-\mu)^T] &=&E[A_0UU^TA_0^T]+E[\varepsilon\varepsilon^T].
\label{eig_A2}
\end{eqnarray}
Note that $E[\varepsilon P_B\varepsilon^T]=\qtilde\sigma^2I_p$ and
$E[\varepsilon \varepsilon^T]=q\sigma^2I_p$. Hence,
$E[(X-\mu)(X-\mu)^T]$ and $E[(X-\mu)P_B(X-\mu)^T]$ have the same
leading $p_0\wedge\ptilde$ eigenvectors as $E[A_0UU^TA_0^T]$ has.
Moreover, we have $\lambda_i=d_i+\qtilde\sigma^2$ and
$\lambda_i^*=d_i+q\sigma^2$, where $d_i$ is the $i^{\rm th}$
eigenvalue of $E[A_0UU^TA_0^T]$. Hence, $\lambda_i^*
-\lambda_i=(q-\qtilde)\sigma^2$ for $i=1,\dots,\ptilde$, which
completes the proof.
\end{proof}

\begin{proof}[\textbf{Proof of Theorem~\ref{thm.tpc}}]
Let ${H_{\ptilde,\qtilde}}(S_n)= (\vec(\widehat A)^T,\vec(\widehat
B)^T)^T$ be the function maps $S_n$ to its tensor principal
components under $(\ptilde,\qtilde)$, which gives
${H_{\ptilde,\qtilde}}(\Sigma)=(\vec(A)^T,\vec(B)^T)^T$ and
$D_{H_{\ptilde,\qtilde}}=\frac{\partial
{H_{\ptilde,\qtilde}}(\Sigma)}{\partial\vec(\Sigma)}$. Note that
$\widehat\Phi(p,q)=\vec(I_m)^T\vec(S_n)$ and
$\Phi(p,q)=\vec(I_m)^T\vec(\Sigma)$. From the weak convergence
$\sqrt{n}(S_n-\Sigma)\stackrel{d}{\rightarrow}N$ and an
application of the delta method, we have, for
$(\ptilde,\qtilde)\leq (p_0,q_0)$,
\begin{eqnarray*}
\sqrt{n}\left(\widehat\Phi(p,q)-\Phi(p,q)\right)&\stackrel{d}{\to}
&\vec(I_m)^T\vec(N),\\
\sqrt{n}\left(\widehat\Phi(\ptilde,\qtilde)-\Phi(\ptilde,\qtilde)\right)
&\stackrel{d}{\to} &D_{\Phi(\ptilde,\qtilde)}\vec(N)
,\\
\sqrt{n}\left({H_{\ptilde,\qtilde}}(S_n)-{H_{\ptilde,\qtilde}}(\Sigma)\right)
&\stackrel{d}{\to}&D_{H_{\ptilde,\qtilde}}\vec(N).
\end{eqnarray*}
The explicit forms of $D_{\Phi(\ptilde,\qtilde)}$ and elements in
$D_{H_{p_0,q_0}}$ are provided in Lemma~\ref{lem.diff}.
\end{proof}

\begin{proof}[\textbf{Proof of Lemma~\ref{lem.diff}}]
For a given pair $(\ptilde,\qtilde)$ with $1\le\ptilde\le p$ and
$1\le\qtilde\le q$, we have, from~(\ref{def2.3_eq1})
and~(\ref{def2.3_eq2}), that $A$ and $B$ satisfy the following
system of stationary equations
\begin{eqnarray*}
\left(\sum_{j=1}^{\qtilde}(b_j\otimes I_p)^T\Sigma(b_j\otimes I_p)\right)a_i&=&\lambda_ia_i, ~i=1,\cdots,\ptilde,\\
\left(\sum_{i=1}^{\ptilde}(I_q\otimes a_i)^T\Sigma(I_q\otimes
a_i)\right)b_j&=&\xi_jb_j, ~j=1,\cdots,\qtilde,
\end{eqnarray*}
where $a_i,b_j,\lambda_i,\xi_j$ depend on
$(\Sigma,\ptilde,\qtilde)$. The indices $i,j$ in the above system of
equations can go beyond $\ptilde$ and $\qtilde$ and up to $p$ and
$q$. But those $a_i$ and $b_j$ with $i>\ptilde$ and $j>\qtilde$ will
not be included in the solution pair $(A,B)$. Note that we have the
following identity, which is due to the definition of $\Phi$ and the
stationary equations:
\begin{equation}
\Phi(\ptilde,\qtilde)=\sum_{i=1}^\ptilde\lambda_i(\Sigma,\ptilde,\qtilde)
=\sum_{j=1}^\qtilde \xi_j(\Sigma,\ptilde,\qtilde).
\end{equation}

We will use the perturbation method (Sibson, Lemma~2.1, 1979; Fine,
1987) to derive the derivatives $D_{\Phi(\ptilde,\qtilde)}$,
$\frac{\partial a_i}{\partial\vec(\Sigma)}$ and $\frac{\partial
b_j}{\partial\vec(\Sigma)}$. Suppose that $\Sigma$ is perturbed to
$\Sigma_\epsilon=\Sigma+\epsilon\dot\Sigma$. Denote the
corresponding system of stationary equations with $\Sigma_\epsilon$
by
\begin{eqnarray}
\left(\sum_{j=1}^{\qtilde}(b_{j,\epsilon}\otimes
I_p)^T\Sigma_\epsilon (b_{j,\epsilon}\otimes I_p)\right)
a_{i,\epsilon}&=&\lambda_{i,\epsilon}~a_{i,\epsilon}, ~i=1,\cdots,\ptilde,\label{eq.26}\\
\left(\sum_{i=1}^{\ptilde}(I_q\otimes
a_{i,\epsilon})^T\Sigma_\epsilon(I_q\otimes
a_{i,\epsilon})\right)b_{j,\epsilon}&=&\xi_j~b_{j,\epsilon},
~j=1,\cdots,\qtilde.\label{eq.27}
\end{eqnarray}
Let their first order expansions be denoted by
\begin{align*}
&\lambda_{i,\epsilon}=\lambda_i+\epsilon \dot\lambda_i+o(\epsilon),
&a_{i,\epsilon}=a_i+\epsilon\dot a_i+o(\epsilon),\\
&\xi_{j,\epsilon}=\xi_j+\epsilon\dot\xi_j+o(\epsilon),
&b_{j,\epsilon}=b_j+\epsilon\dot b_j+o(\epsilon).
\end{align*}
Following the same arguments as in Lemma~2.1 of Sibson (1979) and by
equating the terms involving $\epsilon$ in (\ref{eq.26}) we have,
for $i=1,\cdots,\ptilde$,
\begin{eqnarray}
\dot\lambda_i&=&a_i^T\dot\Sigma_Ba_i,\label{d.a0}\\
\dot a_i &=&\left\{\lambda_iI_p-\sum_{j=1}^{\qtilde}(b_j\otimes
I_p)^T\Sigma(b_j\otimes I_p)\right\}^{+}\dot\Sigma_{B}a_i
\label{d.a},
\end{eqnarray}
where
\begin{eqnarray}
\dot\Sigma_{B}&=&\sum_{j=1}^{\qtilde}\left((\dot b_j\otimes
I_p)^T\Sigma(b_j\otimes I_p)+(b_j\otimes I_p)^T\Sigma(\dot
b_j\otimes I_p)\right)+\sum_{j=1}^{\qtilde}(b_j\otimes
I_p)^T\dot\Sigma(b_j\otimes I_p)\nonumber\\
&=&E[X(\dot BB^T+B\dot B^T)X^T]+\sum_{j=1}^{\qtilde}(b_j\otimes
I_p)^T\dot\Sigma(b_j\otimes I_p) \label{d.sig.b}.
\end{eqnarray}
Since $B_\epsilon^TB_\epsilon=I_\qtilde$, then $\dot B=[~\dot
b_1,\cdots,\dot b_{\qtilde}~]$ must satisfy $\dot B^TB+B^T\dot B=0$.

(a) For $(\ptilde,\qtilde)\le(p_0,q_0)$, the first term of
$\sum_{i=1}^{\ptilde}a_i^T\dot\Sigma_Ba_i$ can be expressed as
\begin{eqnarray}
\sum_{i=1}^{\ptilde}a_i^TE[ X(\dot BB^T+B\dot
B^T)X^T]a_i&=&\sum_{j=1}^{\qtilde}\left(\dot b_j^TE[
X^TP_AX]b_j+b_j^TE[X^TP_AX]\dot b_j\right)\nonumber
\end{eqnarray}
which vanishes by noting that $b_j$ is an eigenvector of
$E[XP_AX^T]$ and $b_j^T\dot b_j=\dot b_j^T b_j=0$. This concludes
that
$\sum_{i=1}^{\ptilde}\dot\lambda_i=\sum_{i=1}^{\ptilde}a_i^T\left(\sum_{j=1}^{\qtilde}(b_j\otimes
I_p)^T\dot\Sigma(b_j\otimes I_p)\right)a_i$ and, hence,
\begin{equation}
D_{\Phi(\ptilde,\qtilde)}=\sum_{i=1}^{\ptilde}\sum_{j=1}^{\qtilde}(b_j\otimes
a_i\otimes b_j\otimes a_i)^T=\vec(P_{B\otimes A})^T.
\end{equation}

(b) Assume now $(\ptilde,\qtilde)=(p_0,q_0)$. To derive the form of
$\frac{\partial a_i}{\partial \vec(\Sigma)}$, we are going to show
that the first term of $\dot\Sigma_{B}$ is zero and conclude
$\dot\Sigma_{B}=\sum_{j=1}^{q_0}(b_j\otimes
I_p)^T\dot\Sigma(b_j\otimes I_p)$. This together with (\ref{d.a})
gives
\begin{eqnarray}
\frac{\partial a_i}{\partial
\vec(\Sigma)}&=&\sum_{j=1}^{q_0}\left(b_j\otimes a_i\otimes
b_j\otimes \{\lambda_iI_p-E[XP_BX^T]\}^{+}\right)^T\nonumber\\
&=&\left\{a_i\otimes \vec(P_{B_0})\otimes
(\lambda_iI_p-E[XP_{B_0}X^T])^{+}\right\}^T(K_{p,q}\otimes I_{pq}).
\end{eqnarray}
as desired, where the second equality follows from
Proposition~\ref{relationship} that $\rmspan(B)=\rmspan(B_0)$ when
$\qtilde=q_0$. To complete the proof, first note that
Proposition~\ref{relationship} ensures the existence of a
nonsingular matrix $\eta$ such that $B_0=B\eta$. From
$X=A_0UB_0^T+\varepsilon$ (remember $\mu=0$) and the independent
structure of $U$ and $\varepsilon$, we can represent the first term
of $\dot\Sigma_B$ as
\begin{eqnarray*}\label{d.sig.1}
E[X(\dot BB^T+B\dot B^T)X^T]&=&E[A_0UB_0^T(\dot BB^T+B\dot
B^T)B_0U^TA_0^T]+\sigma^2\tr(B^T\dot
B+\dot B^TB)I_p\\
&=&E[A_0U\eta^T(B^T\dot B+\dot B^TB)\eta
U^TA_0^T]+\sigma^2\tr(B^T\dot B+\dot B^TB)I_p.
\end{eqnarray*}
The proof is completed by noting that $B^T\dot B+\dot B^TB=0$. The
case of $\frac{\partial b_j}{\partial\vec(\Sigma)}$ can be
established in a similar way.
\end{proof}

\begin{proof}[\textbf{Proof of Corollary~\ref{asymp.proj}}]
Consider the function $F(A,B)=P_{B\otimes A}$ with the corresponding
differential $D_{F(A,B)}$. From Theorem~\ref{thm.tpc} (b) and delta
method, we have
\begin{eqnarray*}
\sqrt{n}~\vec(P_{\widehat B\otimes \widehat A}-P_{B\otimes
A})&=&\sqrt{n}~\vec(F(\widehat A,\widehat B)-F(A,B))\\
&=&D_{F(A,B)}\sqrt{n}\left(\left[\begin{array}{c}\vec(\widehat A)\\
\vec(\widehat B)
\end{array}\right]-\left[\begin{array}{c}\vec(A)\\
\vec(B)
\end{array}\right]\right)+o(\frac{1}{\sqrt{n}})\\
&\stackrel{d}{\rightarrow}&D_{P_{B\otimes A}}\vec(N),
\end{eqnarray*}
where $D_{P_{B\otimes A}}=D_{F(A,B)}D_{H_{\ptilde,\qtilde}}$. When
$(\ptilde,\qtilde)=(p_0,q_0)$, the expression of $D_{P_{B\otimes
A}}$ is obtained by a direct calculation together with
Lemma~\ref{lem.diff} (b) and Theorem~\ref{relationship} (a).
\end{proof}

\begin{proof}[\textbf{Proof of Corollary~\ref{asymp.ratio}}]
Consider the function $F(x,y)=x/y$ with the corresponding
differential $D_{F(x,y)}=(y^{-1},-xy^{-2})^T$. From
Theorem~\ref{thm.tpc} (a) and delta method, we have
\begin{eqnarray*}
\sqrt{n}(\hat\rho(\ptilde,\qtilde)-\rho(\ptilde,\qtilde))&=&\sqrt{n}(F(\widehat\Phi(\ptilde,\qtilde),\widehat\Phi(p,q))
-F(\Phi(\ptilde,\qtilde),\Phi(p,q)))\\
&=&D_{F(\Phi(\ptilde,\qtilde),\Phi(p,q))}
\sqrt{n}\left(\left[\begin{array}{c}\widehat\Phi(\ptilde,\qtilde)\\
\widehat\Phi(p,q)
\end{array}\right]-\left[\begin{array}{c}\Phi(\ptilde,\qtilde)\\
\Phi(p,q)
\end{array}\right]\right)+o(\frac{1}{\sqrt{n}})\\
&\stackrel{d}{\rightarrow}&D_{F(\Phi(\ptilde,\qtilde),\Phi(p,q))}\left[\begin{array}{c}D_{\Phi(\ptilde,\qtilde)}\\
\vec(I_m)^T
\end{array}\right]
\vec(N).
\end{eqnarray*}
A direct calculation gives the expression of the asymptotic variance
$\sigma^2_{\rho(\ptilde,\qtilde)}$.
\end{proof}

\begin{proof}[\textbf{Proof of Theorem~\ref{asymp.test}}]
Under $H_0$, we have from Corollary~\ref{asymp.ratio} that, for $n$
large enough,
\begin{eqnarray*}
P\left(\frac{\sqrt{n}(\hat\rho(\ptilde,\qtilde)-\rho_0)}
{\sigma_{\rho(\ptilde,\qtilde)}}>z_\alpha\right)\leq \alpha.
\end{eqnarray*}
The consistency of $\hat\sigma_{\rho(\ptilde,\qtilde)}^2$ and
Slutsky's theorem complete the proof.
\end{proof}

\begin{proof}[\textbf{Proof of Theorem~\ref{thm.efficiency}}]
Since $(\ptilde,\qtilde)=(p_0,q_0)$, we have $A=A^*$ and $B=B^*$
from Theorem~\ref{thm.2d}, and $\rmspan(A)=\rmspan(A_0)$ and
$\rmspan(B)=\rmspan(B_0)$ from Theorem~\ref{relationship}~(a). Let
$\{a_i:i>p_0\}$ and $\{b_j:j>q_0\}$ be orthogonal bases of
$\rmspan(Q_{A_0})$ and $\rmspan(Q_{B_0})$,
$W_{B,q'}=\sum_{j=1}^{q'}(b_j\otimes I_p\otimes b_j\otimes
I_p)^T$, $q'=1,\cdots,q$ and $W_{A,p'}=\sum_{i=1}^{p'}(I_q\otimes
a_i\otimes I_{q}\otimes a_i)^T$, $p'=1,\cdots,p$. Also define
$M_A=[(a_1\otimes M_{A1}),\cdots,(a_{p_0}\otimes M_{Ap_0})]^T$ and
$M_B=[(b_1\otimes M_{B1}),\cdots,(b_{q_0}\otimes M_{Bq_0})]^T$,
where $M_{Ai}=\{\lambda_iI_p-E[XP_{B_0}X^T]\}^{+}$ and
$M_{Bj}=\{\xi_jI_q-E[X^TP_{A_0}X]\}^{+}$. By using these notations
and from Theorem~\ref{thm.tpc}~(b), we have the limiting
distribution of MPCA
\begin{eqnarray}
\sqrt{n}(\vec(\widehat A,\widehat
B)-\vec(A,B))\stackrel{d}{\rightarrow}M_0W_0\vec(N),\label{asymp.tensor}
\end{eqnarray}
where $M_0=\left[\begin{array}{cc}M_A&0 \\ 0
&M_B\end{array}\right]$, $W_0=\left[\begin{array}{c}W_{B,q_0}
\\ W_{A,p_0}\end{array}\right]$. By Lemma~A.1 below, the limiting distribution
of $({\rm 2D})^2$PCA is derived to be
\begin{eqnarray}
\sqrt{n}(\vec(\widehat A^*,\widehat
B^*)-\vec(A^*,B^*))\stackrel{d}{\rightarrow}M_0(W_0+W_{0+})\vec(N),\label{asymp.2d}
\end{eqnarray}
where $W_{0+}=\left[\begin{array}{c}W_{B,q_0+}
\\  W_{A,p_0+}\end{array}\right]$, $W_{A,p_0+}=W_{A,p}-W_{A,p_0}$, and $W_{B,q_0+}=W_{B,q}-W_{B,q_0}$.

Note that $A=A^*$ and $B=B^*$. To complete the proof, by an
application of delta method, it thus suffices to show
\begin{eqnarray}
{\rm aCov}(\vec(\widehat A^*, \widehat B^*))-{\rm
aCov}(\vec(\widehat A,\widehat B))\geq 0.\label{asymp.efficiency}
\end{eqnarray}
From (\ref{asymp.tensor})-(\ref{asymp.2d}) we are left to show
\begin{eqnarray}
M_0(W_0\Sigma_NW_{0+}^T+W_{0+}\Sigma_NW_0^T+W_{0+}\Sigma_NW_{0+}^T)M_0^T\geq0,
\end{eqnarray}
where $\Sigma_N={\rm Cov}(\vec(N))=(I_{m^2}+K_{m,m})(\Sigma\otimes
\Sigma)$ under normality of $\vec(X)$. We are going to show
$M_0W_0\Sigma_NW_{0+}^TM_0^T=0$. This together with the fact
$M_0W_{0+}\Sigma_NW_{0+}^TM_0^T\geq0$ then establishes the desired
result. Observe that
\begin{eqnarray}
&&M_0W_0\Sigma_NW_{0+}^TM_0^T\nonumber\\
&=&\left[\begin{array}{cc} M_AW_{B,q_0}\Sigma_NW_{B,q_0+}^TM_A^T&
H_AW_{B,q_0}\Sigma_NW_{A,p_0+}^TM_B^T\\
H_BW_{A,p_0}\Sigma_NW_{B,q_0+}^TM_A^T &
H_BW_{A,p_0}\Sigma_NW_{A,p_0+}^TM_B^T\end{array}\right]\nonumber.
\end{eqnarray}
From model~(\ref{model}), $\Sigma=(B_0\otimes
A_0)(T+\sigma^2I_{m_0})(B_0\otimes A_0)^T+\sigma^2Q_{B_0\otimes
A_0}$, where $Q_{B_0\otimes A_0}=Q_{B_0}\otimes
P_{A_0}+P_{B_0}\otimes Q_{A_0}+Q_{B_0}\otimes Q_{A_0}$. This implies
$W_{B,q_0}\Sigma_NW_{B,q_0+}^T=0$ and
$W_{A,p_0}\Sigma_NW_{A,p_0+}^T=0$ and, hence, the diagonal elements
of the above matrix vanish. For the off-diagonal elements, the same
reasoning can be used to deduce that
$(M_AW_{B,q_0})\Sigma_NW_{A,p_0+}^T=0$ and
$(M_BW_{A,p_0})\Sigma_NW_{B,q_0+}^T=0$, which establishes
(\ref{asymp.efficiency}). A direct calculation further gives
$$M_0W_{0+}\Sigma_NW_{0+}^TM_0^T=\sigma^4\left[\begin{array}{cc}(q-q_0)M_A(I_{p^2}+K_{p,p})M_A^T & 0
\\ 0& (p-p_0)M_B(I_{q^2}+K_{q,q})M_B^T\end{array}\right],$$
which equals a zero matrix if and only if $(p_0,q_0)=(p,q)$.
\end{proof}

\noindent{\bf Lemma A.1.} {\it Assume model~(\ref{model}) and
assume that the leading eigenvalues
$\{\lambda_i^*:i=1,\cdots,p_0\}$ and $\{\xi_j^*:j=1,\cdots,q_0\}$
of (2D)$^2$PCA are simple roots. Then, the differentials of
(2D)$^2$PCA components with respect to $\Sigma$ under
$(\ptilde,\qtilde)=(p_0,q_0)$ are given by}
\begin{eqnarray}\label{diff.2d}
\frac{\partial \vec(A^*)}{\partial
\vec(\Sigma)}=M_AW_{B,q}~~~\text{and}~~~ \frac{\partial
\vec(B^*)}{\partial \vec(\Sigma)}=M_BW_{A,p}.
\end{eqnarray}
\begin{proof}
We only derive the differential of $A^*$, where the case of $B^*$ is
similarly obtained. Remember that (2D)$^2$PCA components $A^*$ are
leading eigenvectors of $K_{A^*}=E[XX^T]$ with eigenvalues
$\lambda_i^*$. A standard argument (Sibson, 1979) then gives
\begin{eqnarray}\label{diff.2d.1}
\frac{\partial \vec(a_i^*)}{\partial
\vec(K_{A^*})}&=&a_i^{*T}\otimes\{\lambda_i^*I_p-K_{A^*}\}^{+}\nonumber\\
&=&a_i^T\otimes M_{Ai},
\end{eqnarray}
where the second equality follows from Theorem~\ref{thm.2d} with
$M_{Ai}$ being defined in Theorem~\ref{thm.efficiency}. Turning to
the differential of $K_{A^*}$ with respect to $\Sigma$. It is always
true that
\begin{eqnarray*}
K_{A^*}=\sum_{j=1}^q(b_j\otimes I_p)^T\Sigma(b_j\otimes I_p),
\end{eqnarray*}
where $\{b_j:j>q_0\}$ are defined in the beginning of
Theorem~\ref{thm.efficiency}. Thus, we have
\begin{eqnarray}\label{diff.2d.2}
\frac{\partial \vec(K_{A^*})}{\partial \vec(\Sigma)}&=&W_{B,q}.
\end{eqnarray}
From (\ref{diff.2d.1})-(\ref{diff.2d.2}) and the chain rule, the
proof is completed.
\end{proof}

\newpage

\newpage

\newpage

\newpage

\newpage

\newpage

\end{document}